\documentclass[preprint,12pt,3p]{elsarticle}

\usepackage{amssymb}
\usepackage{amsmath,amsthm}
\usepackage{mathrsfs}
\usepackage{hyperref}
\usepackage{cleveref}
\usepackage[all,cmtip]{xy}
\usepackage{epsf, graphicx}
\usepackage{latexsym,amsfonts,amsbsy,amssymb}
\usepackage{geometry}
\usepackage{titletoc}
\usepackage{color,xcolor}
\usepackage{rotating}
\usepackage{algorithm}
\usepackage{algorithmic}
\usepackage{amscd}

\makeatletter

\@addtoreset{equation}{section} \makeatother
\newtheorem{theorem}{Theorem}[section]
\newtheorem{lemma}[theorem]{Lemma}
\newtheorem{notation}[theorem]{Notation}

\newtheorem{remark}[theorem]{Remark}
\newtheorem{void}[theorem]{}
\newtheorem{definition}[theorem]{Definition}
\newtheorem{proposition}[theorem]{Proposition}

\def\IBr{{\rm IBr}}
\def\Irr{{\rm Irr}}
\def\Ind{{\rm Ind}}
\def\Res{{\rm Res}}

\def\Br{{\rm Br}}
\def\br{{\rm br}}

\def\Aut{{\rm Aut}}

\def\Hom{{\rm Hom}}
\def\det{{\rm det}}
\def\End{{\rm End}}

\def\O{\mathcal{O}}

\def\F{\mathcal{F}}

\def\L{\mathcal{L}}
\def\P{\mathcal{P}}
\def\Z{\mathbb{Z}}

\makeatletter
\def\ps@pprintTitle{%
\let\@oddhead\@empty
\let\@evenhead\@empty
\def\@oddfoot{\reset@font\hfil\thepage\hfil}
\let\@evenfoot\@oddfoot
}
\makeatother

\begin{document}

\begin{frontmatter}

\title{On Morita equivalences with endopermutation source and isotypies}

\author{Xin Huang}


\begin{abstract}
We introduce a new type of equivalence between blocks of finite group algebras called an {\it almost isotypy}. An almost isotypy restricts to a weak isotypy in Brou\'{e}'s original definition \cite[D\'{e}finition 4.6]{Broue1990}, and it is slightly weaker than Linckelmann's version \cite[Definition 9.5.1]{Lin18b}.
We show that a bimodule of two block algebras of finite groups - which has an endopermutation module as a source and which induces a Morita equivalence - gives rise, via slash functors, to an almost isotypy if the character values of a (hence any) source are rational integers. Consequently, if two blocks are Morita equivalent via a bimodule with endopermutation source, then they are almost isotypic. We also explain why the notion of almost isotypies is reasonable.
\end{abstract}

\begin{keyword}
finite groups \sep blocks \sep endopermutation modules \sep Morita equivalences \sep isotypies
\end{keyword}

\end{frontmatter}


\section{Introduction}\label{s1}

Throughout this section $p$ is a prime, $k$ is an algebraically closed field of characteristic $p$, and $\O$ is a complete discrete valuation ring with quotient field $K$ of characteristic $0$ and residue field $k$. Assume that $K$ is a splitting field for all finite groups considered below.

Let $G$ and $H$ be finite groups, $b$ a block of $\O G$ and $c$ a block of $\O H$. Denote by $\Z\Irr_K(G,b)$ the group of generalised characters of $G$ over $K$ associated with the block $b$, and denote by $\Z\IBr_K(G,b)$ the corresponding group of generalised Brauer characters. Following Brou\'{e}, a {\it perfect isometry} between $b$ and $c$ is a group isomorphism 
$\Phi:\Z\Irr_K(H,c)\cong \Z\Irr_K(G,b)$
satisfying certain conditions (see \cite[D\'{e}finition 1.4]{Broue1990} or \cite[Definition 9.2.2]{Lin18b}). By arithmetic properties of a perfect isometry, $\Phi$ induces an isomorphism
$\bar{\Phi}:\Z\IBr_K(H,c)\cong\Z\IBr_K(G,b)$
such that $d_G\circ \Phi=\bar{\Phi}\circ d_H$ (see \cite[Corollary 9.2.7]{Lin18b}). Here $d_G$ and $d_H$ are the usual decomposition maps. Given a $p$-element $u$ of $G$ and a block $e$ of $kC_G(u)$, denote by $\hat{e}$ the unique block of $\O C_G(u)$ that lifts $e$. For any $\chi\in {\rm Cl}_K(G)$ (where ${\rm Cl}_K(G)$ denotes the set of $K$-valued class functions on $G$) associated with the block $b$, define a class function $d_{(G,b)}^{(u,e)}(\chi)$ in ${\rm Cl}_K(C_G(u)_{p'})$ by setting $d_{(G,b)}^{(u,e)}(\chi)(s)=\chi(\hat{e}us)$ for all $p'$-elements $s$ in $C_G(u)$.

\medskip In the rest of this section, we fix the following notation:

\begin{notation}\label{notation:section1 notation}
{\rm	Let $G$ and $H$ be finite groups, $b$ a block of $\O G$ and $c$ a block of $\O H$. Suppose that (up to isomorphism) $b$ and $c$ have a common defect group $P$. Let $i\in (\O Gb)^P$ and $j\in (\O Hc)^P$ be source idempotents. Suppose further that the fusion system of the source algebras $i\O Gi$ and $j\O Hj$ on $P$ are equal. For any subgroup $Q$ of $P$, denote by $e_Q$ the unique block of $kC_G(Q)$ satisfying $\br_Q^{\O G}(i)e_Q\neq 0$ and by $f_Q$ the unique block of $kC_H(Q)$ satisfying $\br_Q^{\O H}(j)f_Q\neq 0$. Denote by $\hat{e}_Q$ and $\hat{f}_Q$ the blocks of $\O C_G(Q)$ and $\O C_H(Q)$ lifting $e_Q$ and $f_Q$, respectively.}
\end{notation}

In \cite{Broue1990}, Brou\'{e} defined the notion of an isotypy. The following definition is Brou\'{e}'s original definition. In order to distinguish different notions, let us call it weak isotypy.

\begin{definition}[cf. {\cite[D\'{e}finition 4.6]{Broue1990}}]\label{defi:isotypy1}
{\rm A {\it weak isotypy} between $b$ and $c$ is a family of perfect isometries
$$\Phi_Q:\Z \Irr_K(C_H(Q),\hat{f}_Q)\cong \Z\Irr_K(C_G(Q),\hat{e}_Q)$$
for every {\bf cyclic} subgroup $Q=\langle u\rangle$ of $P$, such that we have an equality of maps 
$$d_{(G,b)}^{(u,e_Q)}\circ \Phi_1=\bar{\Phi}_{Q}\circ d_{(H,c)}^{(u,f_Q)}$$
from $\Z\Irr_K(H,c)$ to $K\otimes_{\Z}\Z\IBr_K(C_G(Q),\hat{e}_Q)$. We say that the perfect isometry $\Phi_1$ {\it extends to a weak isotypy} between $b$ and $c$, and the family $(\Phi_Q)_{\{1\neq Q({\rm cyclic})\subseteq P\}}$ is called a {\it local system}.
}
\end{definition}

When studying $p$-permutation equivalences, a strong version of isotypy (see {\cite[\S 10]{Lin08}} or \cite[Definition 9.5.1]{Lin18b}) was defined; see also \cite[Definition 15.3]{BP} for an equivalent formulation.  The only known examples of isotypic blocks in the sense of \cite[Definition 9.5.1]{Lin18b} are $p$-permutation equivalent blocks (see  \cite[Theorem 10.1]{Lin08} and \cite[Theorem 1.6 (a)]{BP}) and Galois conjugate blocks (see \cite[Theorem 9.6.1]{Lin18b}).  In order to describe the relationship between Morita equivalences with endopermutation source and isotypies, we define the notion an almost isotypy which is slightly weaker than \cite[Definition 9.5.1]{Lin18b}.  

\begin{definition}\label{defi:isotypy2}
	{\rm  An {\it almost isotypy} between $b$ and $c$ is a family of perfect isometries
		$$\Phi_Q:\Z \Irr_K(C_H(Q),\hat{f}_Q)\cong \Z\Irr_K(C_G(Q),\hat{e}_Q)$$
		for {\bf every} subgroup $Q$ of $P$, with the following properties.

		\begin{enumerate} [\rm (i)] 
			\item (Equivariance) For any isomorphism $\varphi:Q\cong R$ in the common fusion system $\mathcal{F}$ we have ${}^\varphi\Phi_Q=\Phi_R$, where ${}^\varphi\Phi_Q$ is obtained from composing $\Phi_Q$ with the isomorphisms $\Z \Irr_K(C_G(Q),\hat{e}_Q)\cong \Z \Irr_K(C_G(R),\hat{e}_R)$ and $\Z \Irr_K(C_H(Q),\hat{f}_Q)\cong \Z \Irr_K(C_H(R),\hat{f}_R)$ given by conjugation with elements $x\in G$ and $y\in H$ satisfying $\varphi(u)=xux^{-1}=yuy^{-1}$ for all $u\in Q$.
		\item (Compatibility) For any subgroup $Q$ of $P$ and any element $u\in C_P(Q)$, there exists $\varepsilon_{Q,u}\in\{\pm1\}$ such that, setting $R=Q\langle u\rangle$, we have  
		$$d_{(C_G(Q),e_Q)}^{(u,e_R)}\circ \Phi_Q=\varepsilon_{Q,u}\cdot\bar{\Phi}_{R}\circ d_{(C_H(Q),f_Q)}^{(u,f_R)}$$
as maps from $\Z\Irr_K(C_H(Q),\hat{f}_Q)$ to $K\otimes_{\Z}\Z\IBr_K(C_G(R),\hat{e}_R)$. We say that the perfect isometry $\Phi_1$ {\it extends to an almost isotypy} between $b$ and $c$, and the family $(\Phi_Q)_{\{1\neq Q\subseteq P\}}$ is called a {\it local system}.		
			
\end{enumerate}
	}
\end{definition}


The only difference between Definition \ref{defi:isotypy2} and the isotypy in \cite[Definition 9.5.1]{Lin18b} is that we add a sign $\varepsilon_{Q,u}$ in the condition (ii). It is easy to see that from an almost isotypy, one obtains a weak isotypy in Definition \ref{defi:isotypy1} (after possibly replacing $\Phi_Q$ by $-\Phi_Q$ for some non-trivial cyclic subgroups $Q$ of $P$).

An $(\O Gb,\O Hc)$-bimodule $M$ induces, via the tensor functor $(K\otimes_\O M)\otimes_{KH}-$, a $\Z$-linear map $\Phi_M:\Z\Irr_K(H,c)\to \Z\Irr_K(G,b)$. The main result of this paper is as follows:

\begin{theorem}\label{theo:main}
Let $V$ be an indecomposable $\F$-stable endopermutation $\O P$-module, viewed as an $\O\Delta P$-module via the canonical isomorphism $\Delta P\cong P$. Let $M$ be an indecomposable direct summand of the $\O Gb$-$\O Hc$-bimodule 
$$\O Gi\otimes_{\O P}{\rm Ind}_{\Delta P}^{P\times P}(V)\otimes_{\O P}j\O H.$$
Suppose that $M$ induces a Morita equivalence between $\O Gb$ and $\O Hc$. For any $1\neq Q\subseteq P$, let $M_Q$ be a $(\Delta Q, e_Q\otimes f_Q^\circ)$-slashed module attached to $M$ over the group $C_G(Q)\times C_H(Q)$; see Lemma \ref{slash} below. Let $\hat{M}_Q$ be an $\O C_G(Q)\hat{e}_Q$-$\O C_H(Q)\hat{f}_Q$-bimodule with an endopermutation source of determinate 1, such that $k\otimes_\O\hat{M}_Q\cong M_Q$ and that $\hat{M}_Q$ induces a Morita equivalence between $\O C_G(Q)\hat{e}_Q$ and $\O C_H(Q)\hat{f}_Q$; see Proposition \ref{prop:conjugation over O} for the existence of $\hat{M}_Q$. Let $\rho_V:P\to \O$ be the character of $V$. The following hold.
	\begin{enumerate} [\rm (i)] 
		\item If $\Phi_M$ extends to a weak isotypy in Definition \ref{defi:isotypy1}, then the values of $\rho_V$ are rational integers.
		\item If the values of $\rho_V$ are rational integers, then $\Phi_M$ together with the family $(\Phi_{\hat{M}_Q})_{\{1\neq Q\subseteq P\}}$ is an almost isotypy.
		\item Assume that $p\geq 3$ and $P$ is abelian. If the values of $\rho_V$ are rational integers, then there exists a family of signs $(\varepsilon_Q)_{\{1\neq Q\subseteq P\}}$ such that $\Phi_M$ together with the family $(\varepsilon_Q\cdot\Phi_{\hat{M}_Q})_{\{1\neq Q\subseteq P\}}$ is an isotypy in \cite[Defnition 9.5.1]{Lin18b}

		\end{enumerate}
	
\end{theorem}

Statement (i) is due to \cite[Theorem]{HZ} - we will prove it again in this paper to slightly repair the arguments in \cite{HZ}. Statement (ii) generalises \cite[Theorem]{HZ} where we only considered weak isotypies. In the proof of statement (iii) (see \ref{vodi:the proof of iii}), we will give a formula to calculate such a family of signs $(\varepsilon_Q)_{\{1\neq Q\subseteq P\}}$.

If $p=2$, $P=Q_8$ (the quaternion group of order 8) and $k\otimes_\O V$ is a $3$-dimensional endotrivial module, the following proposition shows that Theorem \ref{theo:main} (iii) does not hold. Note that $3$-dimensional endotrivial $kQ_8$-modules exist: for example, the Heller translate of a $5$-dimensional ``exotic" endotrivial module described in  \cite[Proposition 4.3]{Th07}.

\begin{proposition}\label{theo:counterexample}
Let $V$ be an indecomposable $\F$-stable endopermutation $\O P$-module, viewed as an $\O\Delta P$-module via the canonical isomorphism $\Delta P\cong P$. Let $M$ be an indecomposable direct summand of the $\O Gb$-$\O Hc$-bimodule 
$$\O Gi\otimes_{\O P}{\rm Ind}_{\Delta P}^{P\times P}(V)\otimes_{\O P}j\O H.$$
Suppose that $M$ induces a Morita equivalence between $\O Gb$ and $\O Hc$, For any $1\neq Q\subseteq P$, let $M_Q$ be a $(\Delta Q, e_Q\otimes f_Q^\circ)$-slashed module attached to $M$ over the group $C_G(Q)\times C_H(Q)$. Let $\hat{M}_Q$ be any $\O C_G(Q)\hat{e}_Q$-$\O C_H(Q)\hat{f}_Q$-bimodule with an endopermutation source such that $k\otimes_\O\hat{M}_Q\cong M_Q$ and that $\hat{M}_Q$ induces a Morita equivalence between $\O C_G(Q)\hat{e}_Q$ and $\O C_H(Q)\hat{f}_Q$. Assume that $P=Q_8$,
and $V$ is an endotrivial $\O P$-module lifting a $3$-dimensional endotrivial $kP$-module; see \cite{Alp01} for the existence of $V$. Then the values of $\rho_V$ are rational integers, and by Theorem \ref{theo:main} $\Phi_M$ extends to an almost isotypy. But for any family of signs $(\varepsilon_Q)_{\{1\neq Q\subseteq P\}}$, $\Phi_M$ together with the family $(\varepsilon_Q\cdot\Phi_{\hat{M}_Q})_{\{1\neq Q\subseteq P\}}$ is not an isotypy in \cite[Defnition 9.5.1]{Lin18b}.
\end{proposition}

\begin{remark}
	{\rm Keep the notation of Proposition \ref{theo:counterexample}. Although $\Phi_M$ together with $(\varepsilon_Q\cdot\Phi_{\hat{M}_Q})_{\{1\neq Q\subseteq P\}}$ is not an isotypy in the sense of \cite[Defnition 9.5.1]{Lin18b}. We still don't know whether $\Phi_M$ extends to an isotypy in the sense of \cite[Defnition 9.5.1]{Lin18b}, because a local system is not necessarily coming from lifting slashed modules.  However, since $M$ has an endopermutation source, lifting slashed modules is the most natural way to find a potential local system. So Proposition \ref{theo:counterexample} shows that the notion of almost isotypies is reasonable, and it seems to be the optimal solution for this problem.}
\end{remark}

Combined with \cite[Theorem 1.13 (b)]{Kessar_Linckelmann} and \cite[Proposition 7.3.13]{Lin18b}, Theorem \ref{theo:main} immediately implies the following theorem.

\begin{theorem}
Denote by $\bar{b}$ and $\bar{c}$ the image of $b$ and $c$ in $kG$ and $kH$, respectively. Assume that $kG\bar{b}$ and $kH\bar{c}$ are Morita equivalent via a bimodule with endopermutation source. Then $\O Gb$ and $\O Hc$ are almost isotypic in the sense of Definition \ref{defi:isotypy2}. Moreover, if $p\geq 3$ and $P$ is abelian, then $\O Gb$ and $\O Hc$ are isotypic in the sense of \cite[Defnition 9.5.1]{Lin18b}.
\end{theorem}

This paper is more than a proof of Theorem \ref{theo:main}. We take this opportunity to prove many properties of local Morita equivalences induced by slashed modules, such as Propositions \ref{prop:vertex of slashed modules}, \ref{prop:local source idempotents}, \ref{prop:compatibility} and \ref{prop:simple modules bijection}. 
In Section \ref{s2} we review some basic concepts and notation. In Section \ref{section:On character values of endopermutation modules} we prove some auxiliary results on character values of endopermutation modules, and in Section \ref{section:s3 On local Morita equivalences induced by slashed modules}, we prove some properties of local Morita equivalences induced by slashed modules. Then we prove in Section \ref{section4:Lifting local Morita equivalences from k to O} that these local Morita equivalences can be uniquely lifted to $\O$ to satisfy some good properties. In Section \ref{section:On generalised decomposition maps} we compare the generalised decomposation numbers of two blocks which are Morita equivalent via a bimodule with endopermutation source. Theorem \ref{theo:main} and Proposition \ref{theo:counterexample} are proved respectively in Section \ref{section6:proofs} and \ref{section:Proof for Q8}. 

\section{Preliminaries}\label{s2}

Throughout this section, $p$ is a prime, $k$ is an algebraically closed field of characteristic $p$; $\O$ is either a complete discrete valuation ring of characteristic $0$ with residue field $k$, or $\O=k$.
In this section we review some terminology and notation. For an algebra $A$, we denote by $A^{\rm op}$ the opposite algebra of $A$ and by $A^\times$ the group of invertible elements in $A$. For an $A$-module $M$, the symbol $[M]$ denotes the isomorphism class of $M$. Unless otherwise specified, all $\O$-algebras and $\O$-modules considered in this paper are $\O$-free of finite $\O$-rank. For a finite group $G$, we denote by $\Delta G$ the the diagonal subgroup $\{(g,g)|g\in G\}$ of the direct product $G\times G$. Whenever useful, we regard an $\O\Delta G$-module as an $\O G$-module and vice versa via the isomorphism $G\cong \Delta G$ sending $g\in G$ to $(g,g)$.
For finite groups $G$ and $H$, an $(\O G,\O H)$-bimodule $M$ can be regarded as an $\O (G\times H)$-module (and vice versa) via $(g,h)m=gmh^{-1}$, where $g\in G$, $h\in H$ and $m\in M$. If $M$ is indecomposable as an $(\O G,\O H)$-bimodule, then $M$ is indecomposable as an $\O (G\times H)$-module, hence has a vertex (in $G\times H$) and a source. 

Given two $\O$-algebras $A$, $B$ and an $(A,B)$-bimodule $M$, ${\rm End}_\O (M)$ is an $(A\otimes_\O A^{\rm op},B\otimes_\O B^{\rm op})$-bimodule: for any $a_1,a_2\in A$, $b_1,b_2\in B$, $\varphi\in {\rm End}_\O(M)$ and $m\in M$,
$$((a_1\otimes a_2)\cdot\varphi\cdot(b_1\otimes b_2))(m)=a_1\varphi(a_2mb_2)b_1.$$

\begin{void}
{\rm \textbf{The Brauer construction.} Let $G$ be a finite group. We refer to \cite[\S10]{Thevenaz} or \cite[Definition 1.3.1]{Lin18a} for the definition of a $G$-algebra. If $A$ is a $G$-algebra (resp. $\O G$-module), we denote by $A^H$ the $N_G(H)$-subalgebra (resp. $\O N_G(H)$-submodule) of $H$-fixed points of $A$ for any subgroup $H$ of $G$.  For any two $p$-subgroups $Q\leq P$ of $G$, the {\it relative trace map} ${\rm Tr}_Q^P:A^Q\to A^P$, is defined by ${\rm Tr}_Q^P(a)=\sum_{x\in [P/Q]}{}^xa$, where $[P/Q]$ denotes a set of representatives of the left cosets of $Q$ in $P$. We denote by ${\rm Br}_P(A)$ or $A(P)$ the {\it Brauer quotient} of $A$, i.e., the $N_G(P)$-algebra (resp. $kN_G(P)$-module)
\[A^P/(\sum_{Q<P}{\rm Tr}_Q^P(A^Q)+J(\O)A^P),\]
where $J(-)$ denotes the Jacobson radical.
We denote by ${\rm br}_P^A:A^P\to A(P)$ the canonical map, which is called the {\it Brauer homomorphism}. Sometimes we write $\br_Q$ instead of $\br_Q^A$ if no confusion arises.

If $A$ is a $G$-interior algebra (for instance, $A={\rm End}_\O(M)$ for some $\O G$-module $M$), then the Brauer quotient $A(P)$ has a natural structure of $C_G(P)$-interior algebra. For the group algebra $\O G$ (considered as a $G$-interior algebra) and a $p$-subgroup $P$, the Brauer homomorphism can be identified with the $\O$-algebra homomorphism
$(\O G)^P\to kC_G(P)$, $\sum_{g\in G}\alpha_gg\mapsto \sum_{g\in C_G(P)}\bar{\alpha}_gg$, where $\bar{\alpha}_g$ denotes the image of $\alpha_g$ in $k$. If $M$ is an $A$-module then $M$ can be viewed as an $\O G$-module via the structure homomorphism $G\to A^\times$. 
}
\end{void}

\begin{void}\label{weights}

{\rm \textbf{Blocks and Brauer pairs.} Let $G$ be a finite group. By a {\it block} of the group algebra $\O G$, we mean a primitive idempotent $b$ of the center of $\O G$, and $\O Gb$ is called a {\it block algebra} of $\O G$. A {\it defect group} of $b$ is a maximal $p$-subgroup $P$ of $G$ such that ${\rm br}_P^{\O G}(b)\neq 0$.

A {\it Brauer pair} (or {\it subpair}) of the group $G$ is a pair $(Q,e_Q)$, where $Q$ is a $p$-subgroup of $G$ and $e_Q$ is a block of $kC_G(Q)$. The Brauer pair $(Q,e_Q)$ is a {\it$b$-Brauer pair} if $e_Q{\rm br}_Q^{\O G}(b)\neq 0$. For the definition of inclusions of Brauer pairs we refer to \cite[Definition 5.9.7]{Lin18a} or \cite[Definition 6.3.1]{Lin18b}.
Denote by $\hat{e}_Q$ the unique block of $\O C_G(Q)$ lifting $e_Q$.

Let $G$ and $H$ be finite groups. Denote by $-^\circ$ the $\O$-algebra isomorphism $\O H\cong(\O H)^{\rm op}$ sending any $h\in H$ to $h^{-1}$. Let $b$ and $c$ be blocks of $\O G$ and $\O H$ respectively. Clearly $c^\circ$ is a block of $\O H$. Then an $(\O Gb, \O Hc)$-bimodule $M$ can be regarded as an $\O(G\times H)$-module belonging to the block $b\otimes c^\circ$ of $\O(G\times H)$, and vice versa. Here, we identify $b\otimes c^\circ$ to its image under the $\O$-algebra isomorphism $\O G\otimes_\O \O H\cong \O(G\times H)$ sending $g\otimes h$ to $(g,h)$ for any $g\in G$ and $h\in H$.

}
\end{void}

\begin{lemma}\label{lem:Brauer pairs}
	Let $G$ be a finite group, $b$ a block of $\O G$ and $(P,e)$ a maximal $b$-Brauer pair of $G$. For any subgroup $Q$ of $P$ denote by $e_Q$ the unique block of $kC_G(Q)$ such that $(Q,e_Q)\leq (P,e)$. The following hold.
	\begin{enumerate}[{\rm (i)}]
		\item Write $Q'=QC_P(Q)$. Then $(C_P(Q),e_{Q'})$ is an $e_Q$-Brauer pair of $C_G(Q)$.
		\item Let $u\in C_P(Q)$ and $R=Q\langle u\rangle$. Then $(\langle u \rangle,e_R)$ is an $e_Q$-Brauer pair of $C_G(Q)$, and we have $(\langle u\rangle,e_R)\leq (C_P(Q),e_{Q'})$ as $e_Q$-Brauer pairs of $C_G(Q)$. 
	\end{enumerate}
\end{lemma}

\noindent{\it Proof.} Let $i$ be a primitive idempotent in $(\O G)^P$ such that ${\rm br}_P^{\O G}(i)e_P\neq 0$. By the definition of inclusion of Brauer pairs \cite[Definition 5.9.7]{Lin18a} and by \cite[Theorem 5.9.6]{Lin18a}, we have $\br_{Q}^{\O G}(i)e_Q\neq 0$, $\br_{Q'}^{\O G}(i)e_{Q'}\neq 0$ and $\br_R^{\O G}(i)e_R\neq 0$. Note that $C_{C_G(Q)}(C_P(Q))=C_G(Q')$ and $C_{C_G(Q)}(u)=C_G(R)$. Then we see that $\br_{C_P(Q)}^{kC_G(Q)}(e_Q)e_{Q'}=\br_{Q'}^{\O G}(e_Q)e_{Q'}\neq 0$ and $\br_{\langle u\rangle}^{kC_G(Q)}(e_Q)e_R=\br_R^{\O G}(e_Q)e_R\neq 0$. Hence by \cite[Definition 6.3.1]{Lin18b}, both $(C_P(Q),e_{Q'})$ and $(\langle u \rangle,e_Q)$ are $e_Q$-Brauer pairs of $C_G(Q)$.

Consider the Brauer homomorphism $\br_{C_P(Q)}^{kC_G(Q)}: (kC_G(Q))^{C_P(Q)}\to kC_G(Q')$. By standard lifting theorems for idempotents, any primitive idempotent in $kC_G(Q')$ is of the form $\br_{C_P(Q)}^{kC_G(Q)}(j)$ for some primitive idempotent $j\in (kC_G(Q))^{C_P(Q)}$. Hence there exists a primitive idempotent $j\in(kC_G(Q))^{C_P(Q)}$ such that $\br_{Q'}^{\O G}(i)\br_{C_P(Q)}^{kC_G(Q)}(j)\neq 0$ and $\br_{C_P(Q)}^{kC_G(Q)}(j)e_{Q'}\neq 0$. By \cite[Theorem 5.9.6 (ii)]{Lin18a}, we have $\br_{\langle u\rangle}^{kC_G(Q)}(j)e_R\neq 0$. Then by \cite[Definition 5.9.7]{Lin18a}, $(\langle u\rangle,e_R)\leq (C_P(Q),e_{Q'})$.   $\hfill\square$

\begin{void}\label{void:points}

{\rm \textbf{Points of algebras.} Let $G$ be a finite group and $A$ an $\O$-algebra. A {\it point} of $A$ is an $A^\times$-conjugacy class of a primitive idempotent in $A$. Let $I$ be a primitive decomposition of $1_A$ in $A$. The {\it multiplicity} of a point $\alpha$ on $A$ is the cardinal $m_\alpha=|I\cap \alpha|$ of the set $I\cap \alpha$, and it does not depend on the choice of $I$. If $A$ is a $G$-algebra, $Q$ is a $p$-subgroup of $G$, a point $\alpha\in A^Q$ is {\it local} if ${\rm br}_Q^A(\alpha)\neq \{0\}$. In this case, $\alpha$ is called a {\it local point} of $Q$ on $A$. We know that ${\rm br}_Q^A(\alpha)$ is a point of $A(Q)$ and that the correspondence $\alpha\mapsto {\rm br}_Q^A(\alpha)$ induces a bijection between the set of the local points of $Q$ on $A$ and the set of the points of $A(Q)$; see e.g. \cite[Lemma 14.5]{Thevenaz}.

Let $G$ be a finite group and $b$ a block of $\O G$. The map $G\to \O Gb$ sending $g$ to $gb$ induces an interior $G$-algebra structure on $\O Gb$. Let $\alpha$ be a local point of $\langle u\rangle$ on $\mathcal{O}Gb$.
Let $\chi\in {\rm Irr}_K(G, b)$ be a character afforded by an $\O G$-module $M$. We set $\chi(u_\alpha):=\chi(ul)$, where $l\in \alpha$. Note that $\chi({u_\alpha})$ is independent of the choice of $l$; see remarks before \cite[Definition 5.15.2]{Lin18a}.
}
\end{void}

\begin{void}
{\rm
\textbf{Almost source algebras.} For the terminology, notation and basic results of fusion systems, we follow \cite[\S8.1]{Lin18b}. Almost source algebras of a block was introduced by Linckelmann in \cite[Definition 4.3]{Lin08}. Let $G$ be a finite group, $b$ a block of $\O G$ and $P$ a defect group of $b$. An idempotent $i\in (\O Gb)^P$ is called an {\it almost source idempotent} if ${\rm br}_P^{\O G}(i)\neq 0$ for every subgroup $Q$ of $P$, there is a unique block $e_Q$ of $kC_G(Q)$ such that ${\rm br}_Q^{\O G}(i)\in kC_G(Q)e_Q$. The $P$-interior algebra $i\O Gi$ is then called an {\it almost source algebra} of the block $b$. If the almost idempotent $i$ is primitive in $(\O Gb)^P$, then $i$ is called a {\it source idempotent} and $i\O Gi$ is called an {\it source algebra} of $b$ (a classical definition by Puig). By \cite[Proposition 4.1]{Lin08}, there is a canonical Morita equivalence between the block algebra $\O Gb$ and the almost source algebra $i\O Gi$ sending an $\O Gb$-module $M$ to the $i\O Gi$-module $iM$. The choice of an almost source idempotent $i\in (\O Gb)^P$ determines a fusion system $\F$ on $P$ such that for any subgroups $Q$ and $R$ of $P$, the set $\Hom_\F(Q,R)$ is the set of all group homomorphisms $\varphi:Q\to R$ for which there is an element $x\in G$ satisfying $\varphi(u)=xux^{-1}$ for all $u\in Q$ and satisfying $xe_Qx^{-1}=e_{xQx^{-1}}$; see \cite[Remark 4.4]{Lin08} or \cite[\S 8.7]{Lin18b}.  Moreover, a subgroup $Q$ of $P$ is fully $\F$-centralised if and only if $C_P(Q)$ is a defect group of the block $e_Q$ of $kC_G(Q)$. Given a subgroup $Q$ of $P$, by the definition of being fully $\F$-centralised, it is always possible to find a subgroup $R$ of $P$ such that $Q\cong R$ in $\F$ and such that $R$ is fully $\F$-centralised.
}
\end{void}

The following result of Linckelmann explains why we will need to work with fully centralised subgroups and almost source idempotent rather than source idempotents.

\begin{proposition}[{\cite[Proposition 4.5]{Lin08}}]\label{almost}
	Let $G$ be a finite group, $b$ a block of $\O G$, $P$ a defect group of $b$, and $i\in (\O Gb)^P$ an almost source idempotent of $b$ with associated almost source algebra $A=i\O Gi$. Let $\F$ be the fusion system on $P$ determined by $i$. If $Q$ is fully $\F$-centralised subgroup of $P$, then $C_P(Q)$ is a defect group of $kC_G(Q)e_Q$ and ${\rm br}_Q^{\O G}(i)$ is an almost source idempotent of $kC_G(Q)e_Q$ with the associated almost source algebra $A(Q)$. In particular, $kC_G(Q)e_Q$ and $A(Q)$ are Morita equivalent.
\end{proposition}

\begin{lemma}[a variation of {\cite[Proposition 5.1]{Lin08}}]\label{lemma:F-stability of source algebras}
Let $G$ be a finite group, $b$ a block of $\O G$, $P$ a defect group of $b$, and $i\in (\O Gb)^P$ a source idempotent of $b$ with the associated source algebra $A=i\O Gi$. Let $\F$ be the fusion system on $P$ determined by $i$. For any $\varphi\in \Aut_\F(P)$, we have
\begin{enumerate}[{\rm (i)}]
	\item $i\O G\cong {}_\varphi (i\O G)$ as $\O P$-$\O Gb$-bimodules and $\O Gi\cong (\O Gi)_{\varphi}$ as $\O Gb$-$\O P$-bimodules;
	\item  $A\cong {}_\varphi A$ as $\O P$-$A$-bimodules and $A\cong A_{\varphi}$ as $A$-$\O P$-bimodules.
\end{enumerate}
\end{lemma}
\noindent{\it Proof} (a modification of the proof of \cite[Proposition 5.1]{Lin08}). Let $e_P$ be the unique block of $kC_G(P)$ such that $\br_P^{\O G}(i)e_P\neq 0$. Since $\varphi\in \Aut_\F(P)$, there is an element $x\in N_G(P,e_P)$, such that $\varphi(u)=xux^{-1}$ for all $u\in P$. Let $\nu$ be the local point of $P$ on $\O Gb$ containing $i$. By \cite[Definition 8.2.1]{Lin18b}, $P$ is $\F$-centric. Then by \cite[Theorem 8.7.3 (iv)]{Lin18b}, $\nu\cap A$ is the unique local point of $P$ on $A$, and we have $x\nu x^{-1}=\nu$. Since $x^{-1}ix\in \nu$, there is an element $c\in ((\O Gb)^P)^\times$ such that $cic^{-1}=x^{-1}ix$, or equivalently, $xcic^{-1}x^{-1}=i$. The map sending $ia$ to $xcia=ixca$ is a required isomorphism $i\O G\cong {}_\varphi (i\O G)$. The map sending $ai$ to $aic^{-1}x^{-1}=ac^{-1}x^{-1}i$ is a required isomorphism $\O Gi\cong (\O Gi)_{\varphi}$. This proves (i). By the standard Morita equivalence between $\O Gb$ and $A$ we obtain (ii).  $\hfill\square$

\begin{void}\label{section:Fusion-stable endopermutation module}
	
{\rm \textbf{Fusion-stable endopermutation modules.} Let $P$ be a finite $p$-group, $V$ an $\O P$-module. If ${\rm End}_\O(V)$ admits a $P$-stable $\O$-basis under the conjugation action, then $V$ is called an {\it endopermutation} $\O P$-module, as defined in \cite{Dade}.
Let $\F$ be a fusion system on $P$. Let $Q$ be a subgroup of $P$ and $V$ an endopermutation $\O Q$-module. Following \cite[Definition 9.9.1]{Lin18a}, we say that $V$ is {\it $\F$-stable} if for any subgroup $R$ of $Q$ and any morphism $\varphi:R\to Q$ in $\F$, the sets of isomorphism classes of indecomposable direct summands with vertex $R$ of the $kR$-modules ${\rm Res}_R^Q(V)$ and ${}_\varphi V$ are equal (including the possibility that both sets may be empty). This is equivalent to require that ${\rm Res}_R^Q(V)\oplus {}_\varphi V$ is an endopermutation $kR$-module; see \cite[Corollary 6.12]{Dade}.
}
\end{void}

\begin{theorem}[Dade's slashed modules; see e.g. {\cite[Proposition 7.3.7 (i)]{Lin18b}}]\label{theorem:Dade's slahsed modules}
Let $P$ be a finite $p$-group and $V$ an endopermutation $\O P$-module. For any subgroup $Q$ of $P$, and any subgroup $H$ of $P$ satisfying $H\leq N_P(Q)$, there is up to isomorphism a unique endopermutation $kH$-module $V[Q]$ such that $(\End_\O(V))(Q)\cong \End_k(V[Q])$ as $H$-algebras. Moreover, $V[Q]$ has an indecomposable direct summand with vertex $H$.
\end{theorem}

\begin{remark}\label{remark:Defres}
{\rm The $kH$-module $V[Q]$ in Theorem \ref{theorem:Dade's slahsed modules} is called a {\it $Q$-slashed module} attached to $V$ over the group $H$. If $Q$ is a normal subgroup of $P$, a $Q$-slashed module attached to $V$ over the group $P$ is also written as ${\rm Def}_{P/Q}^P(V)$. For this reason, the slash construction is also known as ``deflation--restriction": the $kH$-module $V[Q]$ in Theorem \ref{theorem:Dade's slahsed modules} is isomorphic to ${\rm Def}_{H/Q}^H({\rm Res}_H^P(V))$.}
\end{remark}

\begin{proposition}[{\cite[Proposition 4.1]{Lin15}}]\label{lin154.1}
	Let $A$ and $B$ be almost source algebras of blocks of finite group algebras having a common defect group $P$ and the same fusion system $\mathcal{F}$ on $P$. Let $V$ be an $\mathcal{F}$-stable indecomposable endopermutation $\O P$-module with vertex $P$. Let $M$ be a direct summand the $(A,B)$-bimodule
	$$A\otimes_{\O P}{\rm Ind}_{\Delta P}^{P\times P}(V)\otimes_{\O P}B.$$
	Consider $M$ as an $\O\Delta P$-module via the homomorphism $\Delta P\to A\otimes_\O B^{\rm op}$ sending $(u,u)\in \Delta P$ to $u1_A\otimes u^{-1}1_B$. Then for any non-trivial subgroup $Q$ of $P$, there is a canonical $(A(Q), B(Q))$-bimodule $M_Q$ satisfying ${\rm Br}_{\Delta Q}({\rm End}_\O(M))\cong {\rm End}_k(M_Q)$ as $k$-algebras and as $(A(Q)\otimes_kA(Q)^{\rm op},B(Q)\otimes_kB(Q)^{\rm op})$-bimodules.
\end{proposition}

A Brauer-friendly module is a direct sum of indecomposable modules with compatible fusion-stable endopermutation sources; we refer to \cite[Definition 8]{Biland} for its definition.

\begin{proposition}[see {\cite[Proposition 2.6]{H23}}]\label{brauerfriendly}
	Let $G$ and $H$ be finite groups, and let $b$ and $c$ be blocks of $\O G$ and $\O H$ respectively. Let $M$ be an indecomposable
	$(\O Gb, \O Hc)$-bimodule with an endopermutation source inducing a stable equivalence of Morita type between $\O Gb$ and $\O Hc$. Then $M$ is a Brauer-friendly $\O(G\times H)(b\otimes c^\circ)$-module.
\end{proposition}

Note that although \cite[Proposition 2.6]{H23} is stated over $k$, the proof works over $\O$.

\begin{lemma}[see {\cite[Theorems 18 and 21]{Biland}}]\label{slash}
	Let $G$ be a finite group, $b$ a block of $\O G$, $M$ a Brauer-friendly $\O Gb$-module, $(P,e_P)$ a $b$-subpair, and $H$ a subgroup of $G$ such that $C_G(P)\leq H\leq N_G(P,e_P)$. Then there exists a Brauer-friendly $kHe_P$-module $Sl_{(P,e_P)}^H(M)$ and an isomorphism of $C_G(P)$-interior $H$-algebras
	$$\theta_{(P,e_P)}^H:~{\rm Br}_P({\rm End}_\O(\hat{e}_PM))\cong {\rm End}_k(Sl_{(P,e_P)}^H(M)).$$
\end{lemma}

Following Biland, the pair $(Sl_{(P,e_P)}^H(M),\theta_{(P,e_P)}^H)$ or just the $kHe_P$-module $Sl_{(P,e_P)}^H(M)$ is called a $(P,e_P)$-{\it slashed module} attached to $M$ over the group $H$. This is a generalisation of Dade's slashed modules for endopermutation modules of finite $p$-groups.

\begin{remark}
{\rm In Lemma \ref{slash}, $H$ is assumed to satisfy $C_G(P)\leq H\leq N_G(P,e_P)$, while in \cite[Theorem 18 and 21]{Biland}, $H$ is assumed to satisfy	$PC_G(P)\leq H\leq N_G(P,e_P)$. But  Lemma \ref{slash} follows from the proof of \cite[Theorems 18 and 21]{Biland}.}
\end{remark}

\begin{void}\label{void:vertex subpair}
{\rm  \textbf{Vertex subpairs and sources triples.} Let $G$ be a finite group, $b$ a block of $\O G$, and $M$ an indecomposable $\O Gb$-module. Let $(P,e_P)$ be a $b$-Brauer pair and  $\hat{e}_P$ the unique block of $\O C_G(P)$ lifting $e_P$. If $P$ is a vertex of $M$, and $M$ is isomorphic to a direct summand of the $\O Gb$-module $b\O G\hat{e}_P\otimes_{\O P} V$ for some indecomposable $\O P$-module $V$, then $(P,e_P)$ is called a {\it vertex subpair} of $M$, $V$ is called a {\it source} of $M$ with respect to the vertex subpair $(P,e_P)$, and $(P,e_P,V)$ is called a {\it source triple} of $M$; see \cite[Definition 2]{Biland}.
As in the classical theory of Green, the vertex subpairs and source triples of $M$ exist, and form an orbit under the action of $G$ by conjugation. 
}
\end{void}

It is clear that the notion of source triples is a block theoretic refinement of the concept of vertex-source pairs. The following proposition is a block theoretic refinement of \cite[Theorem 5.6.9]{Lin18a}, and it is also a generalisation of \cite[Proposition 5.3 (b), (c)]{BP}.

\begin{proposition}[a characterisation of vertex subpairs]\label{prop:vertex subpair}
Let $G$ be a finite group, $b$ a block of $\O G$, $M$ an indecomposable $\O Gb$-module, and $(Q,e_Q)$ a $b$-Brauer pair of $G$. Then $(Q,e_Q)$ is contained in a vertex subpair of $M$ if and only if ${\rm Br}_Q(\End_\O (\hat{e}_QM))\neq 0$. 
\end{proposition}

\noindent{\it Proof.} Assume that $(P,e_P)$ is a vertex subpair of $M$ such that $(Q,e_Q)\leq (P,e_P)$. Then by \cite[Lemma 1 (ii)]{Biland}, the $\O P$-module $\hat{e}_PM$ $(=\hat{e}_PbM)$ admits an indecomposable direct summand $V$ with vertex $P$. Since $V$ is indecomposable, there is a primitive idempotent $i\in (\O Gb)^P$ such that $\hat{e}_Pbi=i\hat{e}_Pb=i$ and that $V$ is isomorphic to a direct summand of $iM$. Since $V$ has vertex $P$, by \cite[Theorem 5.6.9]{Lin18a}, we have ${\rm Br}_P(\End_\O (iM))\neq 0$, and hence $\br_P^{\O G}(i)\neq 0$. It follows that $\br_P^{\O G}(i)e_P=\br_P^{\O G}(i\hat{e}_P)=\br_P^{\O G}(ib\hat{e}_P)=\br_P^{\O G}(i)\neq 0$. Since $Q\leq P$, we have ${\rm Br}_Q(\End_\O (iM))\neq 0$. It follows that
$$0\neq{\rm Br}_Q(\End_\O (iM))\cong \Br_{\Delta Q}(i^\circ M^*\otimes_\O iM)\cong \br_{\Delta Q}^{\O(G\times G)}(i^\circ\otimes i)\Br_{\Delta Q}(M^*\otimes_\O M)$$
(at least) as $k$-modules,  where the second isomorphism is by \cite[Lemma 3.9]{Lin08}. By the definition of $(Q,e_Q)\leq (P,e_P)$, we deduce that
 $$\br_{\Delta Q}^{\O(G\times G)}(i^\circ\otimes i)\br_{\Delta Q}^{\O(G\times G)}(\hat{e}_Q^\circ\otimes \hat{e}_Q)=\br_{\Delta Q}^{\O(G\times G)}(i^\circ\otimes i).$$
Hence $${\rm Br}_Q(\End_\O (\hat{e}_QM))\cong \Br_{\Delta Q}(\hat{e}_Q^\circ M^*\otimes_\O \hat{e}_QM)\cong \br_{\Delta Q}^{\O(G\times G)}(\hat{e}_Q^\circ\otimes \hat{e}_Q)\Br_{\Delta Q}(M^*\otimes_\O M)\neq 0,$$ 
proving the ``only if" part.

The proof of the ``~if~" part is inspired by the proof of \cite[Proposition 5.3 (c)]{BP}.  Assume conversely that ${\rm Br}_Q(\End_\O (\hat{e}_QM))\neq 0$. Then the $\O Q$-module $\hat{e}_QM$ admits an indecomposable direct summand $V$ with vertex $Q$.  Since ${\rm Br}_Q(\End_\O (M))\neq 0$, by \cite[Theorem 5.6.9]{Lin18a}, $Q$ is contained in a vertex, say $P$, of $M$. We proceed by induction on the index $|P:Q|$. If $Q=P$, then by \cite[Lemma 1 (ii)]{Biland}, $(P,e_Q)=(Q,e_Q)$ is a vertex subpair of $M$, hence the claim is true. Assume now that $Q<P$. Let $c:={\rm Tr}_{N_G(Q,e_Q)}^{N_G(Q)}(\hat{e}_Q)$.  By \cite[Theorem 6.2.6 (iii)]{Lin18b}, $\hat{e}_Q$ remains a block of $\O N_G(Q,e_Q)$, $c$ is a block of $\O N_G(P)$, and the $\O N_G(Q)c$-$\O N_G(Q,e_Q)\hat{e}_Q$-bimodule $\O N_G(Q)\hat{e}_Q$ and its dual $\hat{e}_Q\O N_G(Q)$ induces a Morita equivalence between $\O N_G(Q)c$ and $\O N_G(Q,e_Q)\hat{e}_Q$. It is well-known that that a splendid Morita equivalence preserves vertices and sources of indecomposable modules; this can be easily proved by using some variation of \cite[Theorem 5.1.16]{Lin18a}.
 Since ${\rm Br}_Q(\End_\O (\hat{e}_QM))\neq 0$, there is an indecomposable direct summand $N$ of the $\O N_G(Q,e_Q)\hat{e}_Q$-module $\hat{e}_QM$ such that $Q$ is contained in a vertex of $N$. By the splendid Morita equivalence above, there is an indecomposable direct summand $N'$ of ${\rm Res}_{N_G(Q)}^{G}(M)$ such that $N=\hat{e}_QN'$ and that a vertex of $N$ is a vertex of $N'$. If $Q$ is a vertex of $N$ (and hence of $N'$), then by the Burry--Carlson--Puig Theorem (\cite[Corollary 5.5.20]{Lin18a}), $Q$ is a vertex of $M$, a contradiction. Thus the $\O N_G(Q,e_Q)\hat{e}_Q$-module $\hat{e}_QN'$ has a vertex $R$ properly containing $Q$. Since $\hat{e}_QN'$ is isomorphic to a direct summand of $\hat{e}_QM$, we have ${\rm Br}_R(\End_\O(\hat{e}_QM))\neq 0$. Hence 
 $$0\neq{\rm Br}_R(\End_\O (\hat{e}_QM))\cong \Br_R(\hat{e}_Q^\circ M^*\otimes_\O \hat{e}_QM)\cong \br_{\Delta R}^{\O(G\times G)}(\hat{e}_Q^\circ\otimes \hat{e}_Q)\Br_R(M^*\otimes_\O M)$$
as $kC_G(R)\br_R(\hat{e}_Q)$-$kC_G(R)\br_R(\hat{e}_R)$-bimodules,  where the second isomorphism is by \cite[Lemma 3.9]{Lin08}. Since $\br_R(\hat{e}_Q)$ is an idempotent in the center of $kC_G(R)$, there exists a block $e_R$ of $kC_G(R)$ such that $e_R\br_R(\hat{e}_Q)\neq 0$ and that 
$$0\neq(e_R^\circ\otimes e_R)\Br_R(M^*\otimes_\O M)\cong \Br_R(\End_\O(\hat{e}_RM)).$$
Note that we have $(Q,e_Q)\leq (R,e_R)$. Applying the induction hypothesis to the $b$-Brauer pair $(R,e_R)$, $(R,e_R)$ (and hence $(Q,e_Q)$) is contained in some vertex subpair of $M$.  $\hfill\square$

\begin{remark}
{\rm In Proposition \ref{prop:vertex subpair}, if $M$ is an indecomposable $p$-permutation $\O Gb$-module, then by \cite[Proposition 5.8.6]{Lin18a} and \cite[Lemma 3.9]{Lin08} we have
	$${\rm Br}_Q(\End_\O (\hat{e}_QM))\cong \End_k(\Br_Q(\hat{e}_QM))\cong {\rm End}_k(e_Q{\rm Br}_Q(M)).$$ So ${\rm Br}_Q(\End_\O (\hat{e}_QM))\neq 0$ if and only if $e_Q{\rm Br}_Q(M)\neq 0$. Hence Proposition \ref{prop:vertex subpair} generalises \cite[Proposition 5.3 (b), (c)]{BP}.
}
\end{remark}


\section{On character values of endopermutation modules}\label{section:On character values of endopermutation modules}

Throughout this section $p$ is a prime, $k$ is an algebraically closed field of characteristic $p$, and $\O$ is a complete discrete valuation ring of characteristic $0$ with residue field $k$. In this section we prove some auxiliary results on character values of endopermutation modules which will be used in the proof of Theorem \ref{theo:main}.

\begin{notation}\label{notation:endopermutation modules}
{\rm Let $P$ be a finite $p$-group and $V$ an endopermutation $\O P$-module such that $\End_\O(V)(P)\neq 0$. Denote by $\rho_V:V\to \O$ the character of $V$ and $\det_V$ the determinant of $V$; see e.g. the paragraph after \cite[Definition 3.1.2]{Lin18a}. Note that $\det_V$ is a linear character of $P$.
For any $u\in P$, by \cite[Proposition 7.3.7 (ii)]{Lin18b}, there is a unique local point $\delta_u$ of $\langle u\rangle$ on $\End_\O(V)$. Denote by $m_u$ the multiplicity of $\delta_u$ on $\End_{\O \langle u\rangle}(V)$. Denote by $\omega_V(u)$ the trace of the $\O$-linear transformation $u\ell\in \End_\O(V)$, where $\ell\in \delta_u$. Clearly $\omega_V(u)$ is independent of the choice of $\ell$. By the proof of \cite[Proposition 52.3]{Thevenaz}, we have $\rho_V(u)=m_u\omega_V(u)$ for any $u\in P$. Also by the proof of \cite[Proposition 52.3]{Thevenaz}, $\rho_V(u)\in \Z$ if and only if $\omega_V(u)=\pm1$. So if $\rho_V(u)\in \Z$, then $\omega_V(u)$ is the sign of $\rho_V(u)$. Assume that $\F$ is a fusion system on $P$ and $V$ is $\F$-stable. It is straightforward to verify that for any $\varphi\in \Hom_\F(\langle u\rangle, P)$, we have $\omega_V(u)=\omega_V(\varphi(u))$.
}
\end{notation}

\begin{proposition}\label{prop:character values and determinant}
Let $P$ be a finite $p$-group and $V$ an endopermutation $\O P$-module having at least one direct summand with vertex $P$. Assume that $p\nmid {\rm rk}_\O(V)$. The following hold.
\begin{enumerate}[{\rm (i)}]
	\item If $\det_V=1$ then the values of $\rho_V$ are in $\Z$. Conversely, if $p\geq 3$ and the values of $\rho_V$ are in $\Z$, then $\det_V=1$.
	\item Assume that $\det_V=1$. For any $u$ and $v$ in $P$ satisfying $\langle u\rangle =\langle v\rangle $, we have $\rho_V(u)=\rho_V(v)$.
\end{enumerate}

\end{proposition}

\noindent{\it Proof.}  (i). Let $W(k)$ be the ring of Witt vectors in $\O$ of $k$ (see e.g. the paragraph before \cite[Definition 1.8]{Kessar_Linckelmann} for the definition). By \cite[Theorem 14.2]{Th07}, there is an endopermutation $W(k)P$-module $V_0$ such that $k\otimes_{W(k)}V_0\cong k\otimes_\O V$. By \cite[Proposition 7.3.12]{Lin18b}, $V_0$ is unique up to tensoring a $W(k)P$-free module of $W(k)$-rank $1$. Since $p\nmid {\rm rk}_\O(V_0)$, by \cite[Corollary 5.3.4]{Lin18a}, we can choose $V_0$ to have determinant $1$ and such a choice is unique up to isomorphism. Then $\O\otimes_{W(k)}V_0$ is an $\O P$-module of determinant $1$ which lifts the $kP$-module $k\otimes_\O V$.

Assume now $\det_V=1$. By \cite[Corollary 5.3.4]{Lin18a} we deduce that $V\cong \O\otimes_{W(k)}V_0$. By \cite[Proposition 7.3.13]{Lin18b}, the values of $\rho_{V_0}$ are in $\Z$. Hence the values of $\rho_V$ are in $\Z$. Conversely, assume that $p\geq 3$ and the values of $\rho_V$ are in $\Z$. Again by \cite[Proposition 7.3.12]{Lin18b}, there is a unique (up to isomorphism) $\O P$-module $T$ of $\O$-rank $1$ such that $T\otimes_\O V\cong \O \otimes_{W(k)}V_0$. Hence we have $\rho_{V_0}=\rho_V\rho_T$ and $\det_{V_0}=(\det_T)^n\det_V$, where $n={\rm rk}_{\O}(V)$. It follows that the values of $\rho_T$ are in $\mathbb{Q}$, and hence in $\Z$. Since $p\geq 3$, the only possible situation is $T\cong\O$. So $\det_V=\det_{V_0}=1$.



(ii). Since $\det_V=1$, by the proof of statement (i), $V$ is defined over $W(k)$. So we may assume that $\O$ is absolutely unramified. Then the statement is proved in the proof of \cite[Proposition 7.3.13]{Lin18b}. $\hfill\square$

\medskip The rest of this section consists of auxiliary results for the proof of Theorem \ref{theo:main} (iii). The reader who does not concern Theorem \ref{theo:main} (iii) can skip the rest of this section. 

\begin{void}
{\rm For a finite $p$-group $P$ and $\Lambda\in \{\O,k\}$, denote by $D_\Lambda(P)$ the Dade group of $P$ over $\Lambda$ (see e.g. \cite[\S3]{Th07}). Let $X$ be a finite $P$-set and let $\Lambda X$ be the corresponding permutation $\Lambda P$-module. Let $\Omega_X(\Lambda)$ be the kernel of the ``augmentation" homomorphism $\Lambda X\to \Lambda$ (mapping every basis element in $X$ to 1). The $\Lambda P$-module $\Omega_X(\Lambda)$ is called a {\it relative syzygy} of $\Lambda$. By a result due to Alperin (see e.g. \cite[Lemma 2.3.3]{Bouc:relative syzygies}), $\Omega_{X}(\Lambda)$ is an endopermutation $\Lambda P$-module.
Let $Q$ be a proper subgroup of $P$, then $P/Q$ (the set of left cosets) is a $P$-set. Let $\Lambda[P/Q]$ be the corresponding permutation $\Lambda P$-module. Since ${\rm Soc}(\Omega_{P/Q}(\Lambda))\subseteq{\rm Soc}(\Lambda[P/Q])\cong \Lambda$, $\Omega_{P/Q}(\Lambda)$ is indecomposable. By \cite[3.2.1]{Bouc:relative syzygies}, $\Omega_{P/Q}(\Lambda)$ has vertex $P$. Let $D^\Omega_\Lambda(P)$ be the subgroup of $D_\Lambda(P)$ generated by all the relative syzygies $\Omega_X(\Lambda)$, where $X$ runs over all non-empty finite $P$-sets. By \cite[Lemma 5.2.3]{Bouc:relative syzygies}, $D^\Omega_\Lambda(P)$ is actually generated by all $\Omega_{P/Q}(\Lambda)$, where $Q$ runs over all proper subgroups of $P$. 
}
\end{void}

\begin{proposition}\label{prop:character values and relative syzygies}
Let $P$ be a finite $p$-group and $V$ an indecomposable endopermutation $\O P$-module with vertex $P$. If $[V]\in  D^\Omega_\O(P)$, then the values of $\rho_V$ are in $\Z$. Conversely, if $p\geq 3$ and the values of $\rho_V$ are in $\Z$, then $[V]\in D_\O^\Omega(P)$. 
\end{proposition}

\noindent{\it Proof.} Assume that $[V]\in D^\Omega_\O(P)$. Since each generator $\Omega_{P/Q}(\O)$ (where $Q$ is a proper subgroup of $P$) is defined over $W(k)$, using \cite[Proposition 7.3.4 (vii)]{Lin18b}, it is easy to see that $V$ is defined over $W(k)$. Hence by \cite[Proposition 7.3.13]{Lin18b}, the values of $\rho_V$ are in $\Z$.   
Assume conversely that $p\geq 3$ and the values of $\rho_V$ in $\Z$. By \cite{Bouc:The Dade group}, we have $D_k(P)=D_k^\Omega(P)$, hence the exists an indecomposable endopermuation $\O P$-module $V'$ such that $[V']\in D_\O^\Omega(P)$ and $k\otimes_\O V'\cong k\otimes_\O V$. By \cite[Proposition 7.3.12]{Lin18b}, there is a unique (up to isomorphism) $\O P$-module $T$ of $\O$-rank $1$ such that $T\otimes_\O V\cong V'$. By the first statement, the values of $\rho_{V'}$ are in $\Z$. It follows that the values of $\rho_T$ are in $\Z$, which in turn implies $T\cong \O$. $\hfill\square$

\begin{lemma}\label{lemma:omege and direct summands}
Let $P$ a finite $p$-group, $W$ an endopermutation $\O P$-module and $V$ an indecomposable direct summand of $W$ with vertex $P$. Then for any $u\in P$ we have 
\begin{enumerate}[{\rm (i)}]
	\item $\omega_V(u)=\omega_W(u)$;
	\item $\rho_V(u)\in \Z$ if and only if $\rho_W(u)\in \Z$.
\end{enumerate}

\end{lemma}

\noindent{\it Proof.}  By \cite[Proposition 7.3.7 (ii)]{Lin18b}, there is a unique local point $\delta_u$ (resp. $\delta'_u$) of $\langle u\rangle$ on $\End_\O(V)$ (resp. $\End_\O(W)$). Let $\ell\in \delta_u$ and $\ell'\in \delta'_u$. Then $\ell(V)$ (resp. $\ell'(W)$) is the unique, up to isomorphism, direct summand of ${\rm Res}_{\langle u\rangle}^P(V)$ (resp. ${\rm Res}_{\langle u\rangle}^P(W)$) with vertex $\langle u\rangle$. Since ${\rm Res}_{\langle u\rangle}^P(V)$ is a direct summand of ${\rm Res}_{\langle u\rangle}^P(W)$, we have $\ell(V)\cong \ell'(W)$. If we regard $u$ as an $\O$-linear transformation of the $\O$-module $\ell(V)$ (resp. $\ell'(W)$), then by definition, $\omega_V(u)$ (resp. $\omega_V(u)$) is the trace of $u\in \End_\O(\ell(V))$ (resp. $u\in \End_\O(\ell'(W))$). So $\omega_V(u)=\omega_W(u)$, whence (i). 

By Notation \ref{notation:endopermutation modules}, $\rho_V(u)\in\Z$ if and only if $\omega_V(u)=\pm1$; $\rho_W(u)\in\Z$ if and only if $\omega_W(u)=\pm1$. Hence (ii) follows by (i).    $\hfill\square$

\begin{lemma}\label{lemma:direct summands and lifting}
Assume that $p\geq 3$. Let $P$ a finite $p$-group, $W$ an endopermutation $kP$-module and $V$ an indecomposable direct summand of $W$ with vertex $P$. Assume that $\hat{W}$ and $\hat{V}$ are endopermutation $\O P$-modules such that $k\otimes_\O \hat{W}\cong W$ and $k\otimes_\O \hat{V}\cong V$. Suppose further that the values of $\rho_{\hat{W}}$ are in $\Z$ and that ${\rm det}_{\hat{V}}=1$. Then $\hat{V}$ is isomorphic to a direct summand of $\hat{W}$.
\end{lemma}

\noindent{\it Proof.}  By \cite[Proposition 7.3.4 (vii)]{Lin18b} and the Krull--Schmidt theorem, there exisits an indecomposable direct summand $\hat{U}$ of $\hat{W}$ such that $k\otimes_\O \hat{U}\cong V$. By Lemma \ref{lemma:omege and direct summands} (ii), the values of $\rho_{\hat{U}}$ are in $\Z$. By Proposition \ref{prop:character values and determinant} (i), this implies $\det_{\hat{U}}=1$.  By the uniquely lifting property with determinate 1 \cite[Corollary 5.3.4]{Lin18a}, we have $\hat{V}\cong \hat{U}$, completing the proof.               $\hfill\square$

\begin{proposition}\label{prop:compativility of omega}
	Let $P$ be a finite $p$-group and $V$ an indecomposable endopermutation $\O P$-module with vertex $P$.  Assume that $p\geq 3$ and the values of $\rho_V$ are in $\Z$. Let $Q_1\cdots,Q_n$ be a sequence of subgroups of $P$ and $u_0=1, u_1, \cdots, u_n$ be a sequence of elements in $P$ such that for any $i\in \{1,\cdots, n\}$,
	$Q_i\leq N_P(\langle u_0,\cdots,u_{i-1}\rangle)$ and $u_i\in Q_i$.
	Let $V_i$ be an indecomposable direct summand of a $\langle u_1,\cdots, u_i\rangle$-slashed module attached to $V$ over $Q_i$ with vertex $Q_i$, and $\hat{V}_i$ an endopermutation $\O Q_i$-module such that $k\otimes_\O\hat{V}_i\cong V_i$ and such that $\det_{\hat{V}_i}=1$. Let $\hat{V}_0=V$. The product $\omega_V(u_1)\omega_{\hat{V}_1}(u_2)\cdots\omega_{\hat{V}_{n-1}}(u_n)$ depends only on the group $\langle u_1,\cdots,u_n\rangle$ and $V$ - it does not depend on a particular choice of the sequence $Q_1,\cdots,Q_n$, $u_1,\cdots, u_n$.
\end{proposition}

\noindent{\it Proof.}  Note that by \cite[Theorem 14.2]{Th07} and  \cite[Corollary 5.3.4]{Lin18a}, $\hat{V}_i$ exists and is unique up to isomorphism. If $n=1$, then the proposition holds by Proposition \ref{prop:character values and determinant} (ii). Since the values of $\rho_V$ are in $\Z$, by Proposition \ref{prop:character values and relative syzygies}, $[V]\in D_{\O}^\Omega(P)$. 

We claim that the statement can be reduced to generators of the group $D_\O^\Omega(P)$.  Let $W$ be another indecomposable endopermutation $\O P$-module with vertex $P$ such that the values of $\rho_W$ are in $\Z$. Let $M$ be an indecomposable direct summand of $V\otimes_\O W$ with vertex $P$ and let $N=V^*$. We can similarly choose $W_i$, $\hat{W}_i$ for $W$, $M_i$, $\hat{M}_i$ for $M$, and $N_i$ and $\hat{N}_i$ for $N$, where $i\in \{1,\cdots,n\}$. Since (by Notation \ref{notation:endopermutation modules}) $\omega_V(u_1)$ is the sign of the integer $\rho_V(u_1)$, we have $\omega_M(u_1)=\omega_{V\otimes_\O W}(u_1)=\omega_V(u_1)\omega_W(u_1)$ and $\omega_{N}(u_1)=\omega_V(u_1)$. Since slash functors are additive (see \cite[Remark 2.8]{H20}) and compatible with tensor products (see \cite[3.2]{H20}), $M_1$ is isomorphic to a direct summand of $V_1\otimes_k W_1$. Note that $\hat{V}_1\otimes_\O \hat{W}_1$ is an endopermutation $\O Q$-module lifting $V_1\otimes_k W_1$ and has determinant $1$. By Proposition \ref{prop:character values and determinant} (i), the values of the character of $\hat{V}_1\otimes_\O \hat{W}_1$ are in $\Z$. Then by Lemma \ref{lemma:direct summands and lifting}, $\hat{M}_1$ is isomorphic to a direct summand of $\hat{V}_1\otimes_\O \hat{W}_1$, and hence $[\hat{M}_1]=[\hat{V}_1\otimes_\O \hat{W}_1]$ in $D_\O(Q)$. So we have $\omega_{\hat{M}_1}(u_2)=\omega_{\hat{V}_1}(u_2)\omega_{\hat{W}_1}(u_2)$. By repeating this process and by using the transitivity of slashed modules (see \cite[Proposition 5.6]{Dade} or \cite[Lemma 22 (i)]{Biland}), we see that $\omega_{\hat{M}_{i-1}}(u_i)=\omega_{\hat{V}_{i-1}}(u_i)\omega_{\hat{W}_{i-1}}(u_i)$ for any $i\in \{1,\cdots, n\}$.
Since slash functors are compatible with duality (see \cite[3.2]{H20}), we have $N_1\cong (V_1)^*$, and hence by the uniquely lifting property with determinate 1 \cite[Corollary 5.3.4]{Lin18a}, $\hat{N}_1\cong(\hat{V}_1)^*$. So we have $\omega_{\hat{N}_1}(s)=\omega_{\hat{V}_1}(s)$ (because the character values $\hat{V}_1$ are in $\Z$). By repeating this process and by using the transitivity of slashed modules (see \cite[Proposition 5.6]{Dade} or \cite[Lemma 22 (i)]{Biland}), we see that $\omega_{\hat{V}_{i-1}}(u_i)=\omega_{\hat{N}_{i-1}}(u_i)$ for any $i\in \{1,\cdots, n\}$. Now we have proved the claim. In other words, to prove the statement, it suffices to consider that case $V=\Omega_{P/R}(\O)$ for some proper subgroup $R$ of $P$.

No we assume that $V=\Omega_{P/R}(\O)$ for some proper subgroup $R$ of $P$. Since $\rho_V+\rho_{\O}$ equals the character of the permutation $\O P$-module $\O[P/R]$, it is straightforward to check that
\begin{equation}\label{equation: omega}
	\omega_V(u_1)=\left\{ \begin{gathered}
		1,~{\rm if}~(P/R)^{\langle u_1\rangle}\neq \emptyset, \\
		-1,~{\rm if}~(P/R)^{\langle u_1\rangle}= \emptyset. \\
	\end{gathered}  \right.
\end{equation}

\textbf{Case 1.} Assume first that $(P/R)^{\langle u_1\rangle}\neq \emptyset$, and hence $\omega_V(u_1)=1$. By the definition of $V_1$, we have $[V_1]=[{\rm Def}_{Q_1/\langle u_1\rangle}^{Q_1}(\Res_{Q_1}^P(V))]$ in $D_k(Q_1)$; see Remark \ref{remark:Defres} for the notation ${\rm Def}_{Q_1/\langle u_1\rangle}^{Q_1}$. By \cite[Corollary 4.1.2 (1) and Lemma 4.2.1 (2)]{Bouc:relative syzygies}, $[V_1]=[\Omega_{(P/R)^{\langle u_1\rangle}}(k)]$ in $D_k(Q_1)$, where we view $(P/R)^{\langle u_1\rangle}$ as a $Q_1$-set. By Lemma \ref{lemma:direct summands and lifting}, we have  $[\hat{V}_1]=[\Omega_{(P/R)^{\langle u_1\rangle}}(\O)]$ in $D_\O(Q_1)$. Hence similar to (\ref{equation: omega}), we have
$$\omega_{\hat{V}_1}(u_2)=\left\{ \begin{gathered}
	1,~{\rm if}~((P/R)^{\langle u_1\rangle})^{\langle u_2\rangle}=(P/R)^{\langle u_1,u_2\rangle}\neq \emptyset, \\
	-1,~{\rm if}~((P/R)^{\langle u_1\rangle})^{\langle u_2\rangle}=(P/R)^{\langle u_1,u_2\rangle}= \emptyset. \\
\end{gathered}  \right.$$

\textbf{Case 2.} Now assume that $(P/R)^{\langle u_1\rangle}= \emptyset$, and hence $\omega_V(u_1)=-1$. By \cite[Corollary 4.1.2 (1) and Lemma 4.2.1 (1)]{Bouc:relative syzygies}, $[V_1]=[k]$ in $D_k(Q_1)$. By Lemma \ref{lemma:direct summands and lifting}, we have  $[\hat{V}_1]=[\O]$ in $D_\O(Q_1)$.  Hence $\omega_{\hat{V}_2}(u_2)=1$. Since $(P/R)^{\langle u_1\rangle}= \emptyset$, we have $(P/R)^{\langle u_1,u_2\rangle}= \emptyset$.

We conclude that in both cases, we have
$$\omega_V(u_1)\omega_{\hat{V}_1}(u_2)=\left\{ \begin{gathered}
	1,~{\rm if}~(P/R)^{\langle u_1,u_2\rangle}\neq \emptyset, \\
	-1,~{\rm if}~(P/R)^{\langle u_1,u_2\rangle}= \emptyset. \\
\end{gathered}  \right.$$ 
By induction and by repeatedly using \cite[Corollary 4.1.2 (1) and Lemma 4.2.1]{Bouc:relative syzygies}, it is routine to show that 
$$\omega_V(u_1)\omega_{\hat{V}_1}(u_2)\cdots\omega_{\hat{V}_{n-1}}(u_n)=\left\{ \begin{gathered}
	1,~{\rm if}~(P/R)^{\langle u_1,\cdots,u_n\rangle}\neq \emptyset, \\
	-1,~{\rm if}~(P/R)^{\langle u_1,\cdots,u_n\rangle}= \emptyset. \\
\end{gathered}  \right.$$ 
This proves the proposition. $\hfill\square$

\begin{proposition}\label{prop:omega and fusion systems}
	Let $P$ be a finite $p$-group and $V$ an indecomposable endopermutation $\O P$-module with vertex $P$.  Assume that $p\geq 3$ and the values of $\rho_V$ are in $\Z$. Assume that $\F$ is a fusion system on $P$ and $V$ is $\F$-stable. Let $Q$ be a subgroup of $P$, $H_Q$ a subgroup of $N_P(Q)$, $V_Q$ an indecomposable direct summand of a $Q$-slashed module attached to $V$ over $H_Q$ with vertex $H_Q$, and $\hat{V}_Q$ an endopermutation $\O H_Q$-module such that $k\otimes_\O\hat{V}_Q\cong V_Q$ and such that $\det_{\hat{V}_Q}=1$. Let $u$ be an element in $H_Q$, $\varphi\in {\rm Hom}_\F(Q\langle u\rangle,P)$, $R=\varphi(Q)$, $v=\varphi(u)$, and $H_R$ a subgroup of $N_P(R)$ containing $v$.  We similarly choose $V_R$ and $\hat{V}_R$. Then $\omega_{\hat{V}_Q}(u)=\omega_{\hat{V}_R}(v)$.
\end{proposition}

\noindent{\it Proof.} Since ${\rm det}_{\hat{V}_Q}=1$, by Proposition \ref{prop:character values and determinant} (i), the values of $\rho_{\hat{V}_Q}$ are in $\Z$. Hence by Notation \ref{notation:endopermutation modules}, $\omega_{\hat{V}_Q}(u)$ is the sign of $\rho_{\hat{V}_Q}(u)$. Again by Notation \ref{notation:endopermutation modules}, $\omega_{\hat{V}_Q}(u)$ is the sign of the character value at $u$ of any indecomposable direct summand of ${\rm Res}_{\langle u\rangle }^{H_Q}(\hat{V}_Q)$ with vertex $\langle u\rangle$. 
Let $W_{Q}$ be an indecomposable direct summand of ${\rm Res}_{Q\langle u\rangle}^{H_Q}(V_Q)$ with vertex $Q\langle u\rangle$ and $\hat{W}_Q$ an endopermutation $\O Q\langle u\rangle$-module such that $k\otimes_\O\hat{W}_Q\cong W_Q$ and such that $\det_{\hat{W}_Q}=1$. Then by Lemma \ref{lemma:direct summands and lifting}, $\hat{W}_Q$ is an indecomposable direct summand of ${\rm Res}_{Q\langle u\rangle }^{H_Q}(\hat{V}_Q)$ with vertex $Q\langle u\rangle$. So $\omega_{\hat{V}_Q}(u)$ is the sign of the character value at $u$ of any indecomposable direct summand of ${\rm Res}_{\langle u\rangle}^{Q\langle u\rangle}(\hat{W}_Q)$ with vertex $\langle u\rangle$. By the definition of $V_Q$, there is a $kH_Q$-module $Y$ such that $V_Q\oplus Y$ is a $Q$-slashed module attached to $V$ over $H_Q$. Using the definition of slashed modules, it is straightforward to check that ${\rm Res}_{Q\langle u\rangle }^{H_Q}(V_Q\oplus Y)$ is a $Q$-slashed module attached to ${\rm Res}_{Q\langle u\rangle}^P(V)$ over $Q\langle u\rangle$. So the isomorphism class of $W_Q$ only depends on the isomorphism class of an indecomposable direct summand of ${\rm Res}_{Q\langle u\rangle}^P(V)$ with vertex $Q\langle u\rangle$.

Similarly, using the discussion above we can choose $W_R$ and $\hat{W}_R$ for $R$ and $v$ instead of $Q$ and $u$, respectively. By the same reason, the isomorphism class of $W_R$ only depends on the isomorphism class of an indecomposable direct summand of ${\rm Res}_{R\langle v\rangle}^P(V)$ with vertex $R\langle v\rangle$. 
Let $S$ be an indecomposable direct summand of ${\rm Res}_{Q\langle u\rangle}^P(V)$ with vertex $Q\langle u\rangle$, and let $T$ be an indecomposable direct summand of ${\rm Res}_{R\langle v\rangle}^P(V)$ with vertex $R\langle v\rangle$. By the definition of the $\F$-stability of $V$, we have $S\cong {}_\varphi T$. From this, we deduce that $W_Q\cong {}_\varphi (W_R)$ as $kQ\langle u\rangle$-modules, and hence $\hat{W}_Q\cong {}_\varphi (\hat{W}_R)$ as $\O Q\langle u\rangle$-modules. This implies $\omega_{\hat{V}_Q}(u)=\omega_{\hat{V}_R}(v)$.

$\hfill\square$

\section{On local Morita equivalences induced by slashed modules}\label{section:s3 On local Morita equivalences induced by slashed modules}
Assume in this section the triple $(K,\O,k)$ is as defined in the beginning of Section \ref{s1}; assume that $G$ and $H$ are finite groups, $b$ is a block of $\O G$ and $c$ is a block of $\O H$; assume that an $(\O Gb,\O Hc)$-bimodule $M$ induces a Morita equivalence between $\O Gb$ and $\O Hc$ and has an endopermutation $\O X$-module $V$.

\begin{notation}\label{notation:3.1}
	
{\rm
 By \cite[Theorem 9.11.2]{Lin18b}, we may assume that 
	\begin{enumerate}[$\bullet$] 
\item the blocks $b$ and $c$ have a common defect group $P$ and a common fusion system $\F$ determined by a source idempotent $i\in (\O Gb)^P$ and also a source idempotent $j\in (\O Hc)^P$; 
\item $X=\Delta P$, and when regarding $V$ as an $\O P$-module via the canonical isomorphism $\Delta P\cong P$, $V$ is $\F$-stable;
\item $M$ is isomorphic to a direct summand of the $(\O Gb,\O Hc)$-bimodule
$$\O Gi\otimes_{\O P}{\rm Ind}_{\Delta P}^{P\times P}(V)\otimes_{\O P}j\O H.$$ 
\end{enumerate}
Moreover, we fix the following notation:
\begin{enumerate}[$\bullet$]
\item  write $A=i\O Gi$ and $B=j\O Hj$;
\item for any subgroup $Q$ of $P$ denote by $e_Q$ the unique block of $kC_G(Q)$ satisfying $\br_Q^{\O G}(i)e_Q\neq 0$ and by $f_Q$ the unique block of $kC_H(Q)$ satisfying $\br_Q^{\O H}(j)f_Q\neq 0$; 
\item denote by $\hat{e}_Q$ and $\hat{f}_Q$ the blocks of $\O C_G(Q)$ and $\O C_H(Q)$ lifting $e_Q$ and $f_Q$, respectively;
\item write $i_Q:=\br_Q^{\O G}(i)$ and $j_Q:=\br_Q^{\O H}(j)$;
\item let $M_Q$ be a $(\Delta Q, e_Q\otimes f_Q^\circ)$-slashed module attached to $M$ over the group $C_G(Q)\times C_H(Q)$. (By Proposition \ref{brauerfriendly}, $M$ is a Brauer-friendly $\O(G\times H)(b\otimes c^{\circ})$-module, hence this makes sense.) Note that by \cite[Lemma 17]{Biland}, $M_Q$ is unique up to isomorphism. 
\end{enumerate}

\begin{theorem}[{\cite[Theorem 1.2 (a)]{H23}}]\label{theorem:H23}
The $(kC_G(Q)e_Q,kC_H(Q)f_Q)$-bimodule $M_Q$ induces a Morita equivalence between $kC_G(Q)e_Q$ and $kC_H(Q)f_Q$ and has an endopermutation module as a source.
\end{theorem}

\begin{proposition}\label{prop:vertex of slashed modules}
Let $Q'=QC_P(Q)$. The pair $(\Delta C_P(Q),e_{Q'}\otimes f_{Q'}^\circ)$ is contained in a vertex subpair of the $kC_G(Q)e_Q$-$kC_H(Q)f_Q$-bimodule $M_Q$, and it is a vertex subpair of $M_Q$ if and only if $Q$ is $\F$-centralised.
\end{proposition}
 
\noindent{\it Proof.} By the definition of a slashed module (see Lemma \ref{slash}), we may regard $M_Q$ as a $k\Delta Q(C_G(Q)\times C_H(Q))$-module on which $\Delta Q$ acts trivially.
By Lemma \ref{lem:Brauer pairs} (i), the pair $(\Delta C_P(Q),e_{Q'}\otimes f_{Q'}^\circ)$ is an $(e_Q\otimes f_Q^\circ)$-Brauer pair of $k(C_G(Q)\times C_H(Q))$. For the first statement, according to Proposition \ref{prop:vertex subpair}, it suffices to show that 
$$\End_k(e_{Q'}M_Qf_{Q'})(\Delta C_P(Q))\neq 0.$$ 
By the transitivity of slash functors \cite[Lemma 22 (i)]{Biland}, we have 
$$\End_k(e_{Q'}M_Qf_{Q'})(\Delta C_P(Q)))\cong \End_\O(\hat{e}_{Q'}M\hat{f}_{Q'})(\Delta Q').$$ 
By \cite[Lemma 2.7]{H23}, $(\Delta P, e_P\otimes f_P^\circ)$ is a vertex subpair of $M$.
Since the $(b\otimes c^\circ)$-Brauer pair $(\Delta Q', e_{Q'}\otimes f_{Q'}^\circ)$ is contained in $(\Delta P,e_P\otimes f_P^\circ)$, again by Proposition \ref{prop:vertex subpair}, we have 
$$\End_\O(\hat{e}_{Q'}M\hat{f}_{Q'})(\Delta Q')\neq 0.$$
This proves the first statement. 


Since $M_Q$ induces a Morita equivalence and has an endopermutation module as a source, by \cite[Theorem 9.11.2 (i)]{Lin18b}, we see that $\Delta C_P(Q)$ is a vertex of $M_Q$ if and only if the order of $C_P(Q)$ equals to the order of a defect group of $kC_G(Q)e_Q$. By \cite[Proposition 8.5.3 (i)]{Lin18b}, this is equivalent to $Q$ being $\F$-centralised. $\hfill\square$
 
\begin{proposition}\label{prop:local source idempotents}
Let $Q'=QC_P(Q)$ and let $(X_Q,e)$ be a vertex subpair of the $kC_G(Q)e_Q$-$kC_H(Q)f_Q$-bimodule $M_Q$ containing $(\Delta C_P(Q),e_{Q'}\otimes f_{Q'}^\circ)$; see Proposition \ref{prop:vertex of slashed modules} for the existence of $(X_Q,e)$. Let $V_Q$ be a source of $M_Q$ with respect to the vertex subpair $(X_Q,e)$. The following hold.
\begin{enumerate}[{\rm (i)}]
	\item Let $X_1$ and $X_2$ be the images of $X_Q$ under the canonical projections $p_1:C_G(Q)\times C_H(Q)\to C_G(Q)$ and $p_2:C_G(Q)\times C_H(Q)\to C_H(Q)$ respectively. Then $p_1$ restricts to an isomorphism $p_1:X_Q\to X_1$, $p_2$ restricts to an isomorphism $p_2:X_Q\to X_2$, and $X_1$ and $X_2$ are defect groups of $kC_G(Q)e_Q$ and $kC_H(Q)f_Q$ respectively. The idempotent $e$ is of the form $b_Q\otimes c_Q^\circ$ where $b_Q$ and $c_Q$ are blocks of $kC_G(QX_1)$ and $kC_H(QX_2)$ respectively.
	\item  There exist source idempotents $s_Q\in (kC_G(Q)e_Q)^{X_1}$ and $t_Q\in (kC_H(Q)f_Q)^{X_2}$ such that $\br_{C_P(Q)}^{kC_G(Q)}(s_Q)e_{Q'}\neq 0$ and $\br_{C_P(Q)}^{kC_H(Q)}(t_Q)f_{Q'}\neq 0$ and that $M_Q$ is isomorphic to a direct summand of 
	$$kC_G(Q)s_Q\otimes_{kX_1}{\rm Ind}_{X_Q}^{X_1\times X_2}(V_Q)\otimes_{kX_2} t_QkC_H(Q).$$ 
	\item Let $\mathcal{F}_Q$ be the fusion system on $X_1$ determined by $s_Q$ and $\mathcal{G}_Q$ be the fusion system on $X_2$ determined by $t_Q$. Write $\psi:X_1\to X_2$ the isomorphism $p_2\circ p_1^{-1}$. Then we have ${}^{\psi}\mathcal{F}_Q=\mathcal{G}_Q$, where ${}^\psi\F_Q$ is the fusion system on $X_2$ induced by $\F_Q$ via the isomorphism $\psi$.
	\item The restriction of $\psi$ to $C_P(Q)$ is the identity map $C_P(Q)\to C_P(Q)$. 
\end{enumerate}

\end{proposition}

\noindent{\it Proof.}  By \cite[Proposition 7.3.10 (i)]{Lin18b}, any source of $M_Q$ has dimension prime to $p$. Then by \cite[Proposition 9.7.1]{Lin18b} we see that $X_1$ is a defect group of $kC_G(Q)e_Q$ and $X_2$ is a defect group of $kC_H(Q)f_Q$. By the definition of a vertex pair, $e$ is a block of the group
$C_{C_G(Q)\times C_H(Q)}(X)=C_{C_G(Q)}(X_1)\times C_{C_H(Q)}(X_2)=C_G(QX_1)\times C_G(QX_2).$ Hence the statement on the idempotent $e$ is now clear, whence (i). Statement (ii) follows from \cite[Lemma 3 (ii)]{Biland}. Starting from statement (ii), using the same argument in the proof of \cite[Theorem 9.11.2 (ii)]{Lin18b}, we obtain statement (iii).  Statement (iv) follows by the choice of $X_Q$ and the definition of $\psi$.   $\hfill\square$

\medskip By \cite[Theorem 9.11.9]{Lin18b}, there is an isomorphism of interior $P$-algebras 
\begin{equation}\label{equation:source algebras}
	f:i\O Gi\cong e(\End_\O(V)\otimes_\O j\O Hj)e
\end{equation}
for some primitive local idempotent $e\in (\End_\O (V)\otimes_\O j\O Hj)^P$. We identify $i\O Gi$ with $e(\End_\O(V)\otimes_\O j\O Hj)e$ via this isomorphism. 
}
\end{notation}
 
\begin{proposition}\label{prop:iMj}
The $(i\O Gi,j\O Hj)$-bimodule $iMj$ is isomorphic to $e(V\otimes_\O j\O Hj)$. Here $e(V\otimes_\O j\O H j)$ is a left $e(\End_\O(V)\otimes_\O j\O Hj)e$-module and we regard it as a left $i\O Gi$-module via the isomorphism $f$; $e(V\otimes_\O j\O H j)$ is also a right $j\O Hj$-module via the multiplication $e(v\otimes b_1)b_2=e(v\otimes b_1b_2)$, for any $v\in V$ and any $b_1,b_2\in j\O Hj$.  
\end{proposition}

\noindent{\it Proof.} The proof of this proposition is to review the choice of $e$ and the construction of $f$ in the proof of \cite[Theorem 9.11.9]{Lin18b}. Since $iMj$ induces a Morita equivalence between $A$ and $B$, we have $A\cong {\rm End}_{B^{\rm op}}(iMj)$ as interior $P$-algebras. Since $iMj$ is isomorphic to a direct summand of 
$$A\otimes_{\O P}\Ind_{\Delta P}^{P\times P}(V)\otimes_{\O P}B,$$ 
we can choose an idempotent $e'$ in $\End_{A\otimes_\O B^{\rm op}}(A\otimes_{\O P}\Ind_{\Delta P}^{P\times P}(V)\otimes_{\O P}B)$ which is a projection to the direct summand $iMj$.
Thus we have an isomorphism
$$f_1:A\cong e'(\End_{B^{\rm op}}(A\otimes_{\O P}\Ind_{\Delta P}^{P\times P}(V)\otimes_{\O P}B))e'$$
of interior $P$-algebras and an isomorphism 
\begin{equation*}\label{equation:iMj}
	iMj\cong e'(A\otimes_{\O P}\Ind_{\Delta P}^{P\times P}(V)\otimes_{\O P}B)
\end{equation*}
of  $A$-modules, where the right side is regarded as an $A$-module via the isomorphism $f_1$.

Since $A$ is a primitive interior $P$-algebra, it follows that the idempotent $e'$ is primitive in $\End_{\O P\otimes_\O B^{\rm op}}(A\otimes_{\O P}\Ind_{\Delta P}^{P\times P}(V)\otimes_{\O P}B)$. So there is an indecomposable summand $W$ of $A$ as an $(\O P,\O P)$-bimodule and an $\End_{\O P\otimes_\O B^{\rm op}}(A\otimes_{\O P}\Ind_{\Delta P}^{P\times P}(V)\otimes_{\O P}B)$-conjugate $e$ of $e'$, such that we have an isomorphism
$$f_2: e'(\End_{B^{\rm op}}(A\otimes_{\O P}\Ind_{\Delta P}^{P\times P}(V)\otimes_{\O P}B))e'\cong e(\End_{B^{\rm op}}(W\otimes_{\O P}\Ind_{\Delta P}^{P\times P}(V)\otimes_{\O P}B))e$$
of interior $P$-algebras. Now we have an isomorphism of $A$-modules
$$iMj\cong e(W\otimes_{\O P}\Ind_{\Delta P}^{P\times P}(V)\otimes_{\O P}B),$$
where the right side is viewed as an $A$-module via the isomorphism $f_2\circ f_1$.  

By \cite[Proposition 2.4.12]{Lin18a}, we have an $(\O P,B)$-bimodule isomorphism ${\rm Ind}_{\Delta P}^{P\times P}(V)\otimes_{\O P} B\cong V\otimes_\O B$. Hence we have an isomorphism 
$$f_3:e(\End_{B^{\rm op}}(W\otimes_{\O P}\Ind_{\Delta P}^{P\times P}(V)\otimes_{\O P}B))e\cong e(\End_{B^{\rm op}}(W\otimes_{\O P} (V\otimes_\O B)))e$$
of interior $P$-algebras.
By \cite[Theorem 8.7.1 (i)]{Lin18b}, $W\cong (\O P)_{\varphi}$ for some $\varphi\in \Aut_\F(P)$. Since $V$ and $B$ are $\F$-stable (see Lemma \ref{lemma:F-stability of source algebras}), we have an isomorphism
$$f_4:e(\End_{B^{\rm op}}(W\otimes_{\O P} (V\otimes_\O B)))e\cong e(\End_{B^{\rm op}}((V\otimes_\O B)))e$$
of interior $P$-algebras.
It is easy to see that we have also an isomorphism 
$$f_5:e(\End_{B^{\rm op}}((V\otimes_\O B)))e\cong e(\End_\O(V)\otimes_\O B)e$$
of interior $P$-algebras.
The isomorphism $f$ is exactly $f_5\circ f_4\circ f_3\circ f_2\circ f_1$. By the above process, we have an isomorphism of $A$-modules $iMj\cong e(V\otimes_\O B)$, where the right side is regarded as an $A$-module via the isomorphism $f$.    $\hfill\square$

\begin{proposition}\label{proposition:i_QM_Qj_Q}
Let $Q$ be a fully $\F$-centralised subgroup of $Q$. Then the $(A(Q), B(Q))$-bimodule $i_QM_Qj_Q$ satisfies the condition that ${\rm Br}_{\Delta Q}({\rm End}_\O(iMj))\cong {\rm End}_k(i_QM_Qj_Q)$ as $k$-algebras and as $(A(Q)\otimes_kA(Q)^{\rm op},B(Q)\otimes_kB(Q)^{\rm op})$-bimodules.
\end{proposition}

\noindent{\it Proof.} By Proposition \ref{almost}, $A(Q)$ and $B(Q)$ are almost source algebras of $kC_G(Q)e_Q$ and $kC_H(Q)f_Q$, respectively. By the definition of a $(\Delta Q,e_Q\otimes f_Q^\circ)$-slashed module, we have ${\rm Br}_Q(\End_\O(M))\cong \End_k(M_Q)$ as interior $(C_G(Q)\times C_H(Q))$-algebras, hence as $(kC_G(Q)e_Q\otimes_k (kC_H(Q)f_Q)^{\rm op}, kC_G(Q)e_Q\otimes_k (kC_H(Q)f_Q)^{\rm op})$-bimodules. The statement follows from the standard Morita equivalences between block algebras and almost source algebras. $\hfill\square$

\begin{lemma}\label{lemma:B(Q)}
For any subgroup $Q$ of $P$ we have ${\rm End}_\O (B)(Q)\cong \End_k(B(Q))$	as $(B(Q),B(Q))$-bimodules and as $k$-algebras.
\end{lemma}

\noindent{\it Proof.} Mimic the proof of \cite[Proposition 5.8.6]{Lin18a}.   $\hfill\square$

\begin{proposition}\label{prop:compatibility}
Let $Q$ be a subgroup of $P$. The following hold.
\begin{enumerate}[{\rm (i)}]
\item There is an isomorphism of interior $C_P(Q)$-algebras 
\begin{equation}\label{equation:source algebras of local blocks}
	\Br_Q(f):A(Q)\cong \br_Q(e)(\End_k(V[Q])\otimes_k B(Q))\br_Q(e)
\end{equation}
for some idempotent $\br_Q(e)\in (\End_k (V[Q])\otimes_k B(Q))^{C_P(Q)}$, where $V[Q]$ is a $Q$-slashed module attached to $V$ over $N_P(Q)$. 
\item Assume that $Q$ is fully $\F$-centralised. The $(A(Q),B(Q))$-bimodule $i_QM_Qj_Q$ is isomorphic to $\br_Q(e)(V[Q]\otimes_k B(Q))$.  Here $\br_Q(e)(V[Q]\otimes_k B(Q))$ is a left $\br_Q(e)(\End_k(V[Q])\otimes_k B(Q))\br_Q(e)$-module and we regard it as a left $A(Q)$-module via the isomorphism $\Br_Q(f)$ in (i); $\br_Q(e)(V[Q]\otimes_k B(Q))$ is also a right $B(Q)$-module via the multiplication $\br_Q(e)(v\otimes b_1)b_2=\br_Q(e)(v\otimes b_1b_2)$, for any $v\in V[Q]$ and any $b_1,b_2\in B(Q)$.  
\end{enumerate}
\end{proposition}

\noindent{\it Proof.} Taking $Q$-Brauer quotient on each side of the interior $P$-algebra isomorphism $f$, see (\ref{equation:source algebras}), we obtain an isomorphism of interior $C_P(Q)$-algebras
$$A(Q)\cong \br_Q^{\End_\O(V)\otimes_\O B}(e)(\End_\O(V)\otimes_\O B)(Q)\br_Q^{\End_\O(V)\otimes_\O B}(e).$$
By \cite[Proposition 5.9.2]{Lin18a} and Theorem \ref{theorem:Dade's slahsed modules}, we have
$$(\End_\O(V)\otimes_\O B)(Q)\cong \End_k(V[Q])\otimes_k B(Q)$$
as interior $C_P(Q)$-algebras, and hence we obtain the desired isomorphism $\Br_Q(f)$ in statement (i).

Next we prove (ii). From the proof of Proposition \ref{lin154.1}, we can easily see that the module ``$M_Q$" in Proposition \ref{lin154.1} is uniquely determined by the module ``$M$", up to isomorphism. So by Proposition \ref{proposition:i_QM_Qj_Q}, it suffices to show that ${\rm Br}_{\Delta Q}({\rm End}_\O(iMj))\cong {\rm End}_k(\br_Q(e)(V[Q]\otimes_k B(Q)))$ as $(A(Q)\otimes_k A(Q)^{\rm op},B(Q)\otimes_kB(Q)^{\rm op})$-bimodules and as $k$-algebras. Indeed, by Proposition \ref{prop:iMj}, we have $iMj\cong e(V\otimes_\O B)$ as $(A,B)$-bimodules. So we have 
\begin{align*}
{\rm Br}_Q({\rm End}_\O(iMj))&\cong \Br_Q(\End_\O(e(V\otimes_\O B)))\\
&\cong \Br_Q(e\End_\O(V\otimes_\O B)e)\\
&\cong \br_Q(e)\Br_Q(\End_\O(V)\otimes_\O \End_\O(B))\br_Q(e)\\
&\cong \br_Q(e)(\End_\O(V)(Q)\otimes_k \End_\O(B)(Q))\br_Q(e)\\
&\cong \br_Q(e)(\End_k(V[Q])\otimes_k \End_k(B(Q)))\br_Q(e)\\
&\cong \br_Q(e)\End_k(V[Q]\otimes_k B(Q))\br_Q(e)\\
&\cong \End_k(\br_Q(e)(V[Q]\otimes_k B(Q)))
\end{align*}
as $(A(Q)\otimes_k A(Q)^{\rm op},B(Q)\otimes_kB(Q)^{\rm op})$-bimodules and as $k$-algebras,
where the third isomorphism is by \cite[Lemma 3.9]{Lin08}, the fourth isomorphism is by \cite[Proposition 5.9.2]{Lin18a}, and the fifth isomorphism is by Theorem \ref{theorem:Dade's slahsed modules} and Lemma \ref{lemma:B(Q)}. Here we use abusively the same notation $\br_Q(e)$ to denote the images of $\br_Q(e)$ in different algebras under the isomorphisms. This completes the proof.   $\hfill\square$

  \begin{void}
 	{\rm For a subgroup $Q$ of $P$, denote by $\L\P(Q,i \O Gi)$ the set of local points of $Q$ on $i\O Gi$, and by $\delta_Q$ the unique local point  of $Q$ on ${\rm End}_\mathcal{O}(V)$. For any $\alpha\in \L\P(Q,i\O Gi)$, let $W_\alpha$ be the simple $A(Q)$-module 
 		$$A(Q){\rm br}_Q^A(l)/J(A(Q)){\rm br}_Q^A(l)$$ for some $l\in \alpha$. Note that $W_\alpha$ is uniquely determined by $\alpha$ up to isomorphism. By \cite[Lemma 14.5]{Thevenaz},   the correspondence $\alpha\mapsto [W_\alpha]$ is a bijection between $\L\P(Q, i\O Gi)$ and the set of isomorphism classes of simple $A(Q)$-modules.

 By \cite[Theorem 7.4.2 (i)]{Lin18b} and the isomorphism (\ref{equation:source algebras}), we obtain a bijection
 		\begin{equation}\label{equation:bijection of local pointed groups}
 			{\cal LP}(Q, i\O Gi)\to {\cal LP}(Q,j\O Hj),\alpha\mapsto \alpha',
 		\end{equation}
 		such that $\alpha$ and $\alpha'$ correspond to each other if an only if both 
 		$$\br_Q(e)(\br_Q^{\End_\O(V)}(\delta_Q)\otimes \br_Q^B(\alpha'))\br_Q(e)$$ 
 		and the image of $\br_Q^A(\alpha)$ under the isomorphism (\ref{equation:source algebras of local blocks}) belong to the same point of the $k$-algebra
 		$$\br_Q(e)(\End_k(V[Q])\otimes_k B(Q))\br_Q(e).$$ 
 		Here we identify the two $C_P(Q)$-interior $N_P(Q)$-algebras $\End_\O(V)(Q)$ and $\End_k(V[Q])$.
 	}
 \end{void}
 
 \begin{proposition}\label{prop:simple modules bijection}
 	Let $Q$ be a fully $\F$-centralised subgroup of $P$. The $(A(Q),B(Q))$-bimodule $i_QM_Qj_Q$ induces a Morita equivalence between $A(Q)$ and $B(Q)$, such that $$i_QM_Qj_Q\otimes_{B(Q)} W_{\alpha'}\cong W_\alpha$$ if $\alpha$ and $\alpha'$ correspond under the bijection (\ref{equation:bijection of local pointed groups}).
 \end{proposition}
 
\noindent{\it Proof.}  By Theorem \ref{theorem:H23}, $M_Q$ induces a Morita equivalence between $kC_G(Q)e_Q$ and $kC_H(Q)f_Q$. By Proposition \ref{almost}, $i_QM_Qj_Q$ induces a Morita equivalence between $A(Q)$ and $B(Q)$. By Proposition \ref{prop:compatibility} (ii), we have $i_QM_Qj_Q\cong \br_Q(e)(V[Q]\otimes_k B(Q))$ as $(A(Q),B(Q))$-bimodules.
Let $l'\in\alpha'$. Then
$$i_QM_Qj_Q\otimes_{B(Q)}W_{\alpha'}\cong \br_Q(e)(V[Q]\otimes_k W_{\alpha'})\cong \br_Q(e)(V[Q]\otimes_k B(Q)\br_Q^B(l')/J(B(Q))\br_Q^B(l')).$$
It is easy to observe that the right side is a simple module of the $k$-algebra 
$$\br_Q(e)(\End_k(V[Q])\otimes_k B(Q))\br_Q(e)$$
corresponding to the point containing $\br_Q(e)(\br_Q^{\End_\O(V)}(\delta_Q)\otimes \br_Q^B(\alpha'))\br_Q(e)$. Now the statement follows from the definition of the bijection (\ref{equation:bijection of local pointed groups}).  $\hfill\square$

  \section{Lifting local Morita equivalences from $k$ to $\O$}\label{section4:Lifting local Morita equivalences from k to O}
Assume in this section the triple $(K,\O,k)$ is as defined in the beginning of Section \ref{s1}.
 For $P$, $Q$ finite $p$-groups, $\mathcal{F}$ a fusion system on $P$ and $\varphi:P\to Q$ a group isomorphism, denote by ${}^\varphi\F$ the fusion system on $Q$ induced by $\F$ via the isomorphism $\varphi$.
 We set $\Delta \varphi:=\{(u,\varphi(u))|u\in P\}$ and whenever useful, we regard an $\O\Delta \varphi$-module $V$ as an $\O P$-module and vice versa via the isomorphism $P\cong \Delta \varphi$ sending $u\in P$ to $(u,\varphi(u))$.

 \begin{lemma}[cf. {\cite[Theorem 1.13]{Kessar_Linckelmann}}]\label{lem:lift-existence}
 	Let $G$ and $H$ be finite groups, $b$ a block of $\O G$ and $c$ a block of $\O H$. Denote by $\bar{b}$ and $\bar{c}$ the images of $b$ and $c$ in $kG$ and $kH$, respectively. For any Morita equivalence (resp. stable equivalence of Morita type) between $kG\bar{b}$ and $kH\bar{c}$ given by an indecomposable bimodule $\bar{M}$ with an endopermutation $kX$-source $\bar{V}$, there is a Morita equivalence (resp. stable equivalence of Morita type) between $\O Gb$ and $\O Hc$ given by an indecomposable bimodule $M$ with an endopermutation $\O X$-source $V$ satisfying the following conditions:
 	\begin{enumerate}[{\rm (i)}]
 		\item $\bar{M}\cong k\otimes_\O M$.
 		\item $\bar{V}\cong k\otimes_\O V$.
 		\item $V$ has determinant $1$. 
 	\end{enumerate}
 \end{lemma}
 \noindent{\it Proof.} The existence of $M$ satisfying conditions (i) and (ii) follows by \cite[Theorem 1.13]{Kessar_Linckelmann}.  By the proof of \cite[Lemma 8.4]{Kessar_Linckelmann}, $M$ also satisfies condition (iii).  $\hfill\square$
 
\begin{remark}\label{remark:on source triples}
{\rm With the notation in Lemma \ref{lem:lift-existence}, if $(X,e,V)$ is a source triple of $M$, then $(X,e,k\otimes_\O V)$ is a source triple of $\bar{M}$; this can be deduced from the definition of a source triple of an indecomposable module (see \ref{void:vertex subpair}). Conversely, since all source triples of $M$ (resp. $\bar{M}$) form a $(G\times H)$-conjugation orbit, if $(Y,f,\bar{W})$ is a source triple of $\bar{M}$, then there exists $(g,h)\in G\times H$ such that $(Y,f,[\bar{W}])={}^{(g,h)}(X,e,[k\otimes_\O V])$, and ${}^{(g,h)}(X,e,V)$ is a source triple of $M$. Write $W:={}^{(g,h)}V$, then $k\otimes_\O W\cong \bar{W}$, and $(Y,f,W)$ is a source triple of $M$.}
\end{remark} 
 
 
 
 The bimodule $M$ in Lemma \ref{lem:lift-existence} is uniquely determined up to isomorphism.
 
 \begin{lemma}\label{lem:lift-uniqueness}
 	Let $G$ and $H$ be finite groups, $b$ a block of $\O G$ and $c$ a block of $\O H$. Let $M$ and $N$ be indecomposable $\O Gb$-$\O Hc$-bimodules with endopermutation source and satisfying $k\otimes_\O M\cong k\otimes_\O N$. Assume that both $M$ and $N$ induce stable equivalences of Morita type between $\O Gb$ and $\O Hc$ and both $M$ and $N$ have sources of determinant $1$. Then $M\cong N$.  
 \end{lemma}
 
 \noindent{\it Proof.} Denote by $\bar{b}$ and $\bar{c}$ the images of $b$ and $c$ in $kG$ and $kH$, respectively. Since a source of $M$ (resp. $N$) is of determinant $1$, any source of $M$ (resp. $N$) is of determinant $1$. Let $P$ be a defect group of $b$ and $Q$ a defect group of $c$. Let $i\in (\O Gb)^P$ and $j\in (\O Hc)^Q$ be source idempotents. Denote by $\F$ the fusion system on $P$ determined by $i$. By \cite[Theorem 9.11.2]{Lin18b}, there are isomorphisms $\varphi:P\to Q$ and $\psi:P\to Q$ and indecomposable endopermutation $\O P$-modules $V$ and $W$, such that $M$ is isomorphic to a direct summand of the $\O Gb$-$\O Hc$-bimodule
 $$\O Gi\otimes_{\O P}\Ind_{\Delta \varphi}^{P\times Q}(V)\otimes_{\O Q}j\O H,$$
 and such that $N$ is isomorphic to a direct summand of the $\O Gb$-$\O Hc$-bimodule
 $$\O Gi\otimes_{\O P}\Ind_{\Delta \psi}^{P\times Q}(W)\otimes_{\O Q}j\O H.$$
 Since $k\otimes_\O M\cong k\otimes_\O N$, by \cite[Lemma 2.7]{BKL}, we see that $\varphi^{-1}\circ\psi\in\Aut_\F(P)$ and that $k\otimes_\O V\cong k\otimes_\O W$ as $kP$-modules. Since both $V$ and $W$ are of determinant $1$, by \cite[Lemma 28.1 (b)]{Thevenaz}, $V\cong W$ as $\O P$-modules. Since $N^*$ is isomorphic to a direct summand of
 $$\O Hj\otimes_{\O Q}\Ind_{\Delta \psi^{-1}}^{Q\times P}({}_{\psi^{-1}}(W^*))\otimes_{\O P} i\O G$$
 (see \cite[page 81]{BKL}), by \cite[Lemma 2.6]{BKL}, $M\otimes_{\O Hb} N^*$ is isomorphic to a direct summand of
 $$\O Gi\otimes_{\O P}\Ind_{\Delta \psi^{-1}\circ \varphi}^{P\times P} (U)\otimes_{\O P} i\O G,$$
 where $U$ is an indecomposable direct summand with vertex $P$ of 
 $$V\otimes_\O  {}_{\varphi\circ\psi^{-1}} (W^*)\cong V\otimes_\O{}_{\varphi\circ\psi^{-1}} (V^*).$$
 By \cite[Theorem 9.11.2 (iii)]{Lin18b}, $V$ is $\F$-stable, hence we have $V^*\cong {}_{\varphi\circ\psi^{-1}} (V^*)$ and $U\cong\O$. So $M\otimes_{\O Hc} N^*$ is a $p$-permutation $\O(G\times H)$module which lifts the $p$-permutation $k(G\times H)$-module $(k\otimes_\O M)\otimes_{kH\bar{c}}(k\otimes_\O N^*)\cong kG\bar{b}\oplus \bar{S}$, where $\bar{S}$ is a projective $kG\bar{b}$-$kH\bar{c}$-bimodule. By \cite[Theorem 5.10.2 (iv)]{Lin18a}, we have $M\otimes_{\O Hc} N^*\cong \O Gb\oplus S$, where $S$ is a projective $\O Gb$-$\O Hc$-bimodule lifting $\bar{S}$. Now it is easy to see that $M\cong N$. $\hfill\square$

 \medskip 

 Assume in the rest of this section $G$ and $H$ are finite groups, $b$ is a block of $\O G$ and $c$ is a block of $\O H$; assume that an $(\O Gb,\O Hc)$-bimodule $M$ induces a Morita equivalence between $\O Gb$ and $\O Hc$ and has an endopermutation $\O X$-module $V$. So we can continue to use Notation \ref{notation:3.1}. 
 
For a finite group $G$ and any element $g\in G$, we use abusively the same notation $c_g$ to denote various group homomorphisms induced by $g$-conjugation. For example, if $g$-conjugation induces an isomorphism between two subgroups $Q$ and $R$ of $G$, then $g$-conjugation also induces an isomorphism between $C_G(Q)$ and $C_G(R)$; both these isomorphisms will be denoted by $c_g$.
 
 \begin{proposition}\label{prop:conjugation over k}
 	For any isomorphism $c_g=c_h:Q\cong R$ in the fusion system $\F$, where $g\in G$ and $h\in H$, we have ${}_{c_g^{-1}} (M_Q)_{c_h^{-1}}\cong M_R$ as $kC_G(R)e_R$-$kC_H(R)f_R$-bimodules.
 \end{proposition}
 
 \noindent{\it Proof.} By the transitivity of slash functors \cite[Lemma 22 (ii)]{Biland}, ${}_{c_g^{-1}} (M_Q)_{c_h^{-1}}$ is another $(\Delta R,e_R\otimes f_R^\circ)$-slashed module attached to $M$ over $C_G(R)\times C_H(R)$. Then the statement follows by \cite[Lemma 17 (ii)]{Biland}. $\hfill\square$

 \begin{proposition}\label{prop:conjugation over O}
 	For any non-trivial subgroup $Q$ of $P$, there exists an $\O C_G(Q)\hat{e}_Q$-$\O C_H(Q)\hat{f}_Q$-bimodule $\hat{M}_Q$ with an endopermutation source of determinate 1, such that $k\otimes_\O\hat{M}_Q\cong M_Q$ and that $\hat{M}_Q$ induces a Morita equivalence between $\O C_G(Q)\hat{e}_Q$ and $\O C_H(Q)\hat{f}_Q$. Such a bimodule $\hat{M}_Q$ is unique up to isomorphism. Moreover, the following hold:
 	\begin{enumerate}[{\rm (i)}]
 		\item $\hat{M}_Q$ has a vertex subpair $(X_Q,e)$ containing the $(e_Q\otimes f_Q^\circ)$-Brauer pair $(\Delta C_P(Q),e_{QC_P(Q)}\otimes f_{QC_P(Q)}^\circ)$ and an endopermutation source $\hat{V}_Q$ of determinant $1$ with respect to the vertex subpair $(X_Q,e)$.
 		\item for any isomorphism $c_g=c_h:Q\cong R$ in the fusion system $\F$, where $g\in G$ and $h\in H$, we have ${}_{c_g^{-1}} (\hat{M}_Q)_{c_h^{-1}}\cong \hat{M}_R$ as $\O C_G(R)\hat{e}_R$-$\O C_H(R)\hat{f}_R$-bimodules.
 	\end{enumerate}
 	
 \end{proposition}
 
 \noindent{\it Proof.} The existence of $\hat{M}_Q$ is according to Lemma \ref{lem:lift-existence}, and the uniqueness of $\hat{M}_Q$ follows by Lemma \ref{lem:lift-uniqueness}. Since a source of $\hat{M}_Q$ has determinate 1, any source of $\hat{M}_Q$ has determinant 1. By Proposition \ref{prop:vertex of slashed modules} and Remark \ref{remark:on source triples},  statement (i) holds. By Proposition \ref{prop:conjugation over k} and Lemma \ref{lem:lift-uniqueness}, statement (ii) holds.   $\hfill\square$
 
 
  \section{On generalised decomposition numbers}\label{section:On generalised decomposition maps}
 
Assume in this section the triple $(K,\O,k)$ is as defined in the beginning of Section \ref{s1}; assume that $G$ and $H$ are finite groups, $b$ is a block of $\O G$ and $c$ is a block of $\O H$; assume that an $(\O Gb,\O Hc)$-bimodule $M$ induces a Morita equivalence between $\O Gb$ and $\O Hc$ and has an endopermutation $\O X$-module $V$. So we can continue to use the notation in Notation \ref{notation:3.1}.
 
 \begin{void}\label{void:simple modules of local blocks}
 	{\rm  Let $Q$ be a subgroup of $P$. A local point $\alpha$ of $Q$ on $\mathcal{O}Gb$ is said to be {\it associated to the Brauer pair $(Q,e_Q)$}, if
 		$e_Q{\rm br}^{\mathcal{O}G}_Q(\alpha)\neq 0$. Denote by $\mathcal{LP}(Q,e_Q)$ the set of local points of $Q$ on $\O Gb$ associated to $(Q,e_Q)$.
 		If $Q$ is a cyclic group generated by $u$, then we also write $(Q,e_Q)$ as $(u,e_u)$ and called it a {\it Brauer element}. Let $S_\alpha$ be the simple $kC_G(Q)e_Q$-module 
 		$$kC_G(Q){\rm br}_Q^{\O G}(l)/(J(kC_G(Q)){\rm br}_Q^{\O G}(l))$$
 		for some $l\in \alpha$. By \cite[Lemma 14.5]{Thevenaz}, the correspondence $\alpha\mapsto [S_\alpha]$ is a bijection between $\L\P(Q, e_Q)$ and the set of isomorphism classes of simple $k C_G(Q)e_Q$-modules. Denote by $\varphi_\alpha$ the Brauer character afforded by the simple module $S_\alpha$. So we have 
 		$$\IBr_K(C_G(Q),e_Q)=\{\varphi_\alpha~|~\alpha\in \mathcal{LP}(Q,e_Q)\}.$$
 	}
 	\end{void}
 	
 	 \begin{lemma}\label{lemma:bijection between pointed groups on blocks and source algebras}
 		Let $Q$ be a fully $\F$-centralised subgroup of $P$. Then there is a bijection between $\L\P(Q,e_Q)$ and the set of local points of $Q$ in $i\O Gi$, sending $\alpha\in \L\P(Q,e_Q)$ to $\alpha\cap i\O Gi$. Similarly, there is a bijection between $\L\P(Q,f_Q)$ and the set of local points of $Q$ in $j\O Hj$, sending $\alpha'\in \L\P(Q,e_Q)$ to $\alpha’\cap j\O Hj$.
 	\end{lemma}
 	
 	\noindent{\it Proof.} By \cite[Proposition 8.7.3 (ii)]{Lin18b}, $\alpha\cap i\O Gi\neq \emptyset$. Now the statement follows from \cite[Proposition 4.12]{Thevenaz}.  $\hfill\square$
 	
 	\begin{proposition}\label{prop:bijections between simple modules of local blocks}
 		Keep the notation of \ref{void:simple modules of local blocks}.	Let $Q$ be a fully $\F$-centralised subgroup of $P$. Then $M_Q$ induces a Morita equivalence between $kC_G(Q)e_Q$ and $kC_H(Q)f_Q$, such that 
 		$$M_Q\otimes_{kC_H(Q)f_Q} S_{\alpha'}\cong S_\alpha$$
 		if $\alpha\cap i\O Gi$ and $\alpha'\cap j\O Hj$ correspond under the bijection (\ref{equation:bijection of local pointed groups}).
 	\end{proposition}
 	
 	\noindent{\it Proof.} This follows from Proposition \ref{prop:simple modules bijection} and the standard Morita equivalences between block algebras and almost source algebras. $\hfill\square$

 	\begin{void}\label{void:decomposition numbers}
 	{\rm Let $\chi\in \Irr_K(G,b)$,  and let $N$ be a simple $KG$-module affording $\chi$. By \cite[Corollary 4.4]{Pointed} (or \cite[Theorem 5.15.3]{Lin18a}), for any $u\in P$ we have
 		\[d_{(G,b)}^{(u,e_u)}(\chi ) = \sum\limits_{\alpha  \in \mathcal{LP}(u,e_u)} {\chi(u_\alpha){\varphi _\alpha }},\] where the notation $\chi(u_\alpha$) is defined in \ref{void:points}. Let $N'$ be a simple $KHc$-module corresponding to $N$ under the Morita equivalence induced by $K\otimes_\O M$ and denote by $\chi'$ the character afforded by $N'$.  Hence $\Phi_M(\chi')=\chi$ (recall that $\Phi_M$ is the $\Z$-linear map $\Z\Irr_K(H,c)\to \Z\Irr_K(G,b)$ induced by the $(K Gb,K Hc)$-bimodule $K\otimes_\O M$ via the tensor product over $KH$). Let $(u, f_u)$ be a Brauer element contained in $(P, f_P)$. Similarly, we have
 		\[d_{(H,c)}^{(u,f_u)}(\chi') = \sum\limits_{\alpha' \in \mathcal{LP}(u,f_u)} {\chi(u_{\alpha'}){\varphi_{\alpha'}}}.\]
 		
 	Assume that $\langle u\rangle$ is fully $\F$-centralised. Let $\alpha\in \L\P(u,e_u)$ and $\alpha'\in \L\P(u,f_u)$ such that $\alpha\cap i\O Gi$ and $\alpha\cap j\O Hj$ correspond under the bijection (\ref{equation:bijection of local pointed groups}). By \cite[Theorem 7.4.3]{Lin18b} we have \[\chi(u_\alpha)=\omega_V(u)\chi'(u_{\alpha'});\]
 	see Notation \ref{notation:endopermutation modules} for the meaning of $\omega_V(u)$. The above formula describes the relationship of the generalised decomposition numbers/matrices of the two blocks $b$ and $c$.
 	}
 \end{void}
 
\begin{remark}
{\rm Since sometimes we regard the $\O \Delta P$-source $V$ of $M$ as an $\O P$-module via the canonical isomorphism $P\cong\Delta P$, we used the notation ``$\omega_V(u)$" in \ref{void:decomposition numbers}. Sometimes if we consider $V$ as an $\O \Delta P$-module, we will write ``$\omega_V((u,u))$". However, we have $\omega_V(u)=\omega_{V}((u,u))$. } 
\end{remark}

 \section{Proof of Theorem \ref{theo:main}}\label{section6:proofs}
 
Let $(K,\O,k)$ be as defined in the beginning of Section \ref{s1}. We first prove an auxiliary lemma.
 
 
 \begin{lemma}\label{lem:F-centralised only}
Keep the notation of Notation \ref{notation:section1 notation}. Assume that we have a family of perfect isometries
$$\Phi_Q:\Z \Irr_K(C_H(Q),\hat{f}_Q)\cong \Z\Irr_K(C_G(Q),\hat{e}_Q)$$
for every cyclic subgroup $Q$ of $P$, with the following properties:
\begin{enumerate} [\rm (i)] 
\item  for any cyclic subgroup $Q$ of $P$ and any isomorphism $\varphi:Q\cong R$ in the fusion system $\mathcal{F}$, we have ${}^\varphi\Phi_Q=\Phi_R$;
\item  for $v\in P$ such that $\langle v\rangle$ is fully $\F$-centralised, we have  
$d_{(G,b)}^{(v,e_v)}\circ \Phi_1=\varepsilon_{1,v}\cdot\bar{\Phi}_{\langle v\rangle}\circ d_{(H,c)}^{(v,f_v)}$ for some $\varepsilon_{1,v}\in\{\pm1\}$.
\end{enumerate}	
For any $u\in P$, there exists an isomorphism $\varphi:\langle u\rangle\to R$ in $\F$ such that $R$ is fully $\F$-centralised. Write $v=\varphi(u)$.
Then we have $d_{(G,b)}^{(u,e_u)}\circ \Phi_1=\varepsilon_{1,v}\cdot\bar{\Phi}_{\langle u\rangle}\circ d_{(H,c)}^{(u,f_u)}$. 
 \end{lemma}
 
 	\noindent{\it Proof.} By condition (i), for any cyclic subgroup $Q$ of $P$ and any isomorphism $\varphi:Q\cong R$ in $\mathcal{F}$, we have ${}^\varphi\bar{\Phi}_Q=\bar{\Phi}_R$, where ${}^\varphi\bar{\Phi}_Q$ is obtained from composing $\bar{\Phi}_Q$ with the isomorphisms $\Z \IBr_K(C_G(Q),\hat{e}_Q)\cong \Z \IBr_K(C_G(R),\hat{e}_R)$ and $\Z \IBr_K(C_H(Q),\hat{f}_Q)\cong \Z \IBr_K(C_H(R),\hat{f}_R)$ given by conjugation with elements $x\in G$ and $y\in H$ satisfying $\varphi(u)=xux^{-1}=yuy^{-1}$ for all $u\in Q$. This is because $\bar{\Phi}_Q$ and $\bar{\Phi}_R$ are uniquely determined by $\Phi_Q$ and $\Phi_R$, respectively.
 	

Write $Q:=\langle u\rangle$. Assume that $x\in G$ and $y\in H$ satisfying $\varphi(w)=xwx^{-1}=ywy^{-1}$ for all $w\in Q$. The equality ${}^\varphi\bar{\Phi}_Q=\bar{\Phi}_R$ means that $\bar{\Phi}_Q$ sends a class function $\lambda' \in K\otimes_\Z\Z\IBr_K(C_H(Q),\hat{f}_Q)$ to $\lambda\in K\otimes_\Z\Z\IBr_K(C_G(Q),\hat{e}_Q)$ if and only if $\bar{\Phi}_R$ sends a class function $(\lambda')^y\in K\otimes_\Z \Z\IBr_K(C_H(R),\hat{f}_R)$ to $\lambda^x\in K\otimes_\Z\Z\IBr_K(C_G(R),\hat{e}_R)$. Here $\lambda^x$ denotes the class function sending an element $g\in C_G(R)_{p'}$ to $\lambda(x^{-1}gx)\in K$, and $(\lambda')^y$ is defined similarly. 
 	
For any $\chi'\in \Z\Irr_K(H,c)$ and $\chi:=\Phi_1(\chi')$, since $R=\langle v\rangle$ is fully $\F$-centralised, by the assumption (ii), the class function $d_{(G,b)}^{(v,e_v)}(\chi)$ equals the class function
 $\varepsilon_{1,v}\cdot\bar{\Phi}_R(d_{(H,c)}^{(v,f_v)}(\chi'))$. In other words, $\bar{\Phi}_R$ sends $d_{(H,c)}^{(v,f_v)}(\chi')$ to $ \varepsilon_{1,v}\cdot d_{(G,b)}^{(v,e_v)}(\chi)$.
A straightforward calculation shows that if we denote the class function $d_{(H,c)}^{(u,f_u)}(\chi'): C_H(Q)_{p'}\to K$ by $\lambda'$, then $(\lambda')^y$ is the class function $d_{(H,c)}^{(v,f_v)}(\chi'):C_H(R)_{p'}\to K$.
 Similarly, if we denote the class function $d_{(G,b)}^{(u,e_u)}(\chi)$ by $\lambda$, then $\lambda^x=d_{(G,b)}^{(v,e_v)}(\chi)$. 
 Since $\bar{\Phi}_R((\lambda')^y)=\varepsilon_{1,v}\cdot\lambda^x$, by the previous paragraph, we have $\bar{\Phi}_Q(\lambda')=\varepsilon_{1,v}\cdot\lambda$. In other words, we have $d_{(G,b)}^{(u,e_u)}(\Phi_1(\chi'))=\varepsilon_{1,v}\cdot\bar{\Phi}_{\langle u\rangle}(d_{(H,c)}^{(u,f_u)}(\chi'))$.    $\hfill\square$
 
\medskip Assume in the rest of this section we are in the context of Theorem \ref{theo:main}. So we can continue to use the notation in Notation \ref{notation:3.1} and Scetion \ref{section:On generalised decomposition maps}. 



\begin{void}
{\rm \textbf{Proof of Theorem \ref{theo:main} (i).}
The key points of the proof are contained in \cite[2.11, 2.12 and 3.5]{HZ}. We prove it again to slightly repair the arguments there. Let $u$ be an element of $P$ such that $\langle u\rangle$ is fully $\F$-centralised. We define a $K$-linear map
$$I_{p'}^u:K\otimes_\Z\Z\IBr_K(C_H(u),f_u)\to K\otimes_\Z \Z\IBr_K(C_G(u),e_u)$$
which sends $\varphi_{\alpha'}$ to $\omega_V(u)\varphi_\alpha$ if the local points $\alpha\cap i\O Gi$ and $\alpha'\cap j\O Hj$ correspond under the bijection (\ref{equation:bijection of local pointed groups}). Then using the three equalities in \ref{void:decomposition numbers}, it is straightforward to check that
\begin{equation}\label{equation:generalised decomposition maps1}
d_{(G,b)}^{(u,e_u)}(\Phi_M(\chi'))=I_{p'}^u(d_{{(H,c)}}^{(u,f_u)}(\chi'))
\end{equation}
for any $\chi'\in \Irr_K(C_H(u),f_u)$.
Note that the map $I_{p'}^u$ is uniquely determined by $\Phi_M$, because $d_{(H,c)}^{(u,f_u)}$ extends linearly to a surjective map $K\otimes_\Z\Z\Irr_K(H,c)\to K\otimes_\Z \Z\IBr_K(C_G(u),e_u)$.

Assume that the perfect isometry $\Phi_M$ extends to a weak isotypy in Definition \ref{defi:isotypy1}.
Then there is a perfect isometry
$$\Phi_{\langle u\rangle}: \Z\Irr_K(C_{H}(u),\hat{f}_{\langle u\rangle})\to \Z\Irr_K(C_{G}(u),\hat{e}_{\langle u\rangle})$$
such that 
$$d_{(G,b)}^{(u,e_u)}(\Phi_M(\chi'))=\bar{\Phi}_{\langle u\rangle}(d_{{(H,c)}}^{(u,f_u)}(\chi')).$$
By the uniqueness of $I_{p'}^u$, we have $\bar{\Phi}_{\langle u\rangle}(\varphi_{\alpha'})=I_{p'}^u(\varphi_{\alpha'})=\omega_V(u)\varphi_\alpha$ for any $\alpha'\in \L\P(u,f_u)$ and $\alpha\in \L\P(u,e_u)$ such that $\alpha'\cap j\O Hj$ and $\alpha\cap i\O Gi$ correspond under the bijection (\ref{equation:bijection of local pointed groups}).
Note that $\bar{\Phi}_{\langle u\rangle}:\Z\IBr_K(C_H(u),\hat{f}_{\langle u\rangle})\to \Z\IBr_K(C_G(u),\hat{e}_{\langle u\rangle})$ is a $\Z$-linear isomorphism. This forces $\omega_V(u)=\pm1$. 

For any $u\in P$, there exists an isomorphism $\varphi$ in $\F$ such that $\langle \varphi(u)\rangle$ is fully $\F$-centralised. Write $v=\varphi(u)$. Denote by $T_u$ (resp. $T_v$)
a direct summand of ${\rm Res}_{\langle u\rangle}^P(V)$ (resp. ${\rm Res}_{\langle v\rangle}^P(V)$) with vertex $\langle u\rangle$ (resp. $\langle v\rangle$).
By the previous paragraph, the character value of $T_v$ at $v$ is $1$. Since $V$ is $\F$-stable, we have $T_u\cong {}_{\varphi^{-1}} (T_v)$, which in turn implies that the character value of $T_u$ at $u$ is $\pm1$. Hence $\omega_V(u)=\pm1$ for all $u\in P$. Then by Notation \ref{notation:endopermutation modules}, we see that $\rho_V(u)\in \Z$ for all $u\in P$.   $\hfill\square$
}
\end{void}

\begin{void}\label{void:The proof of (ii) implying (iv)}
	{\rm  \textbf{Proof of Theorem \ref{theo:main} (ii).}  For any non-trivial subgroup $Q$ of $P$, we set $\hat{M}_Q$ and $\hat{V}_Q$ be as in Proposition \ref{prop:conjugation over O}. For $Q=1$, set $\hat{M}_Q:=M$ and $\hat{V}_Q:=V$. Since $\hat{M}_Q$ induces a Morita equivalence between $\O C_Q(Q)\hat{e}_Q$ and $\O C_H(Q)\hat{f}_Q$, the $\Z$-linear isomorphism 
		$$\Phi_{\hat{M}_Q}:\Z\Irr_K(C_H(Q),\hat{f}_Q)\to \Z\Irr_K(C_G(Q),\hat{e}_Q)$$
		is a perfect isometry. Since $k\otimes_\O\hat{M}_Q\cong M_Q$, the resulting group isomorphism
		$$\bar{\Phi}_{\hat{M}_Q}:\Z\IBr_K(C_H(Q),\hat{f}_Q)\to \Z\IBr_K(C_G(Q),\hat{e}_Q)$$
		is exactly induced by the tensor functor $M_Q\otimes_{kC_H(Q)}-$.
		Consider the family $(\Phi_{\hat{M}_Q})_{Q\subseteq P}$ of perfect isometries. By Proposition \ref{prop:conjugation over O} (iii), the equivariance condition (i) in Definition \ref{defi:isotypy2} holds. For any $u\in P$ such that $\langle u\rangle$ is fully $\F$-centralised, by Proposition \ref{prop:bijections between simple modules of local blocks},
		$\bar{\Phi}_{\hat{M}_{\langle u\rangle}}(\varphi_{\alpha'})=\varphi_\alpha$ for any pair $(\alpha,\alpha')\in \mathcal{LP}(u,f_u)\times\mathcal{LP}(u,f_u)$ such that $\alpha\cap i\O Gi$ and $\alpha'\cap j\O Hj$ correspond under the bijection (\ref{equation:bijection of local pointed groups}).  
		Then using the three equalities in \ref{void:decomposition numbers}, it is straightforward to check that
		\begin{equation*}
			d_{(G,b)}^{(u,e_u)}(\Phi_M(\chi'))=\omega_V(u)\bar{\Phi}_{\hat{M}_{\langle u\rangle}}(d_{{(H,c)}}^{(u,f_u)}(\chi'))
		\end{equation*}
		for any $\chi'\in \Irr_K(kC_H(u),f_u)$. Hence for any ``fully $\F$-centralised element" $u$, we have $d_{(G,b)}^{(u,e_u)}\circ\Phi_M=\omega_V(u)\bar{\Phi}_{\hat{M}_{\langle u\rangle}}\circ d_{{(H,c)}}^{(u,f_u)}$. By Lemma \ref{lem:F-centralised only} and the last sentence in Notation \ref{notation:endopermutation modules}, we have 
		\begin{equation}\label{equation:generalised decomposition maps}
			d_{(G,b)}^{(u,e_u)}\circ\Phi_M=\omega_V(u)\bar{\Phi}_{\hat{M}_{\langle u\rangle}}\circ d_{{(H,c)}}^{(u,f_u)}
		\end{equation}
		for any $u\in P$. As we mentioned before, if we view $V$ as an $\O \Delta P$-module, we write $\omega_V((u,u))$ instead of $\omega_V(u)$. However, we have $\omega_V(u)=\omega_{V}((u,u))$.
	
	Let $Q$ be a subgroup of $P$,  $u$ an element in $C_P(Q)$, and $R$ the group $Q\langle u\rangle$. Consider the block $\hat{e}_Q\otimes \hat{f}_Q^\circ$ of $\O(C_G(Q)\times C_H(Q))$. By Lemma \ref{lem:Brauer pairs} (ii), $(\Delta \langle u\rangle, e_R\otimes f_R^\circ)$ is an $(\hat{e}_Q\otimes \hat{f}_Q^\circ)$-Brauer pair. Consider the $\O C_G(Q)\hat{e}_Q$-$\O C_H(Q)\hat{f}_Q$-bimodule $\hat{M}_Q$. By the transitivity of slash functors \cite[Lemma 22 (i)]{Biland}, $M_R$ is a $(\Delta \langle u\rangle, e_R\otimes f_R^\circ)$-slashed module attached to $\hat{M}_Q$ over $C_G(R)\times C_H(R)$. By Proposition \ref{prop:conjugation over O} (i) and Proposition \ref{prop:character values and determinant} (i), the character values of a source of $\hat{M}_Q$ are rational integers. Recall from Proposition \ref{prop:conjugation over O} that $\hat{M}_Q$ has an $\O X_Q$-source $\hat{V}_Q$ and $X_Q$ contains $\Delta C_P(Q)$. By Proposition \ref{prop:local source idempotents}, the bimodule $\hat{M}_Q$ and the bimodule $M$ satisfy some similar conditions. Thus we can apply the equality (\ref{equation:generalised decomposition maps}) to the $\O C_G(Q)\hat{e}_Q$-$\O C_H(Q)\hat{f}_Q$-bimodule $\hat{M}_Q$ instead of the $\O Gb$-$\O Hc$-bimodule $M$, and we obtain  
	\begin{equation}\label{equation:compatibility R}
	d_{(C_G(Q),e_Q)}^{(u,e_R)}\circ \Phi_{\hat{M}_Q}=\omega_{\hat{V}_Q}((u,u))\bar{\Phi}_{\hat{M}_R}\circ d_{(C_H(Q),f_Q)}^{(u,f_R)}.
\end{equation}

 To prove that the compatibility condition (ii) in Definition \ref{defi:isotypy2} holds, it suffices to show that $\omega_{\hat{V}_Q}((u,u))=\pm1$. By the assumption, the values of $\rho_{\hat{V}_1}:=\rho_{V}$ are in $\Z$.  Hence by Notation \ref{notation:endopermutation modules}, $\omega_{\hat{V}_1}(u,u)=\omega_V(u,u)=\pm 1$ for any $u\in P$. For any non-trivial subgroup $Q$ of $P$, by our choice of $\hat{V}_Q$ and by Proposition \ref{prop:character values and determinant} (i), the values of $\rho_{\hat{V}_Q}$ are in $\Z$. Again by Notation \ref{notation:endopermutation modules}, $\omega_{\hat{V}_Q}(u,u)=\pm 1$ for any $u\in C_P(Q)$.
	 Summarising the above, the family $(\Phi_{\hat{M}_Q})_{Q\subseteq P}$ of perfect isometries together with the family of signs $(\varepsilon_{Q,u}:=\omega_{\hat{V}_Q}(u,u))_{Q\subseteq P;u\in C_P(Q)}$ forms an almost isotypy.   $\hfill\square$
}
\end{void}

\begin{void}\label{vodi:the proof of iii}
	{\rm  \textbf{Proof of Theorem \ref{theo:main} (iii).} We continue to use the notation in \ref{void:The proof of (ii) implying (iv)}. Since $P$ is abelian, any subgroup $Q$ of $P$ is fully $\F$-centralised. Hence both $kC_G(Q)e_Q$ and $kC_H(Q)f_Q$ have $P$ as a defect group. Then by Proposition \ref{prop:vertex of slashed modules}, $\Delta P$ is a vertex of $M_Q$. Let $V_Q:=k\otimes_\O \hat{V}_Q$. Recall that the notation $\hat{V}_Q$ is from Proposition \ref{prop:conjugation over O}, which is a source of $\hat{M}_Q$, and hence $V_Q$ is a source of $M_Q$. By Proposition \ref{prop:conjugation over O} and Proposition \ref{prop:vertex of slashed modules}, $(\Delta P, e_P\otimes f_P^\circ, V_Q)$ is a source triple of $M_Q$. For any subgroup $Q$ of $P$, we can choose a sequence of elements $u_0=1,u_1,\cdots,u_n$ in $P$ such that
		$Q=\langle u_1,\cdots,u_n\rangle$. For $i\in \{1,\cdots,n\}$, let $V_i$ be an indecomposable direct summand of a $\Delta\langle u_1,\cdots, u_i\rangle$-slashed module attached to $V$ over $\Delta C_P(\langle u_0,\cdots,u_{i-1}\rangle)=\Delta N_P(\langle u_1,\cdots, u_i\rangle)=\Delta P$ with vertex $\Delta P$.  By \cite[Lemma 3 (iii)]{Bilandadv}, $V_{\langle u_1,\cdots, u_i\rangle}$ is isomorphic to a direct summand of $V_i$.  Let $\hat{V}_i$ be an endopermutation $\O \Delta P$-module such that $k\otimes_\O\hat{V}_i\cong V_i$ and such that $\det_{\hat{V}_i}=1$. By Lemma \ref{lemma:direct summands and lifting}, $\hat{V}_{\langle u_1,\cdots, u_i\rangle}$ is isomorphic to a direct summand of $\hat{V}_i$. Hence by Lemma \ref{lemma:omege and direct summands}, we have 
		\begin{equation}\label{equation:aaa}
		\omega_{\hat{V}_{\langle u_1,\cdots, u_i\rangle}}=\omega_{\hat{V}_i}.
		\end{equation}
		Let $\varepsilon_Q:=\omega_V((u_1,u_1))\omega_{\hat{V}_1}((u_2,u_2))\cdots\omega_{\hat{V}_{n-1}}((u_n,u_n))$. By Proposition \ref{prop:compativility of omega}, $\varepsilon_Q$ depends only on $Q$ and $V$. 
		
		Consider the family $(\varepsilon_Q\Phi_{\hat{M}_Q})_{Q\subseteq P}$ of perfect isometries. By \ref{void:The proof of (ii) implying (iv)}, the equivariance condition (i) in Definition \ref{defi:isotypy2} holds for the family $(\Phi_{\hat{M}_Q})_{Q\subseteq P}$. Hence by Proposition \ref{prop:omega and fusion systems}, the equivariance condition (i) also holds for the family $(\varepsilon_Q\Phi_{\hat{M}_Q})_{Q\subseteq P}$. For any subgroup $Q$ of $P$ and $u$ an element in $C_P(Q)$, writing $R:=Q\langle u\rangle$, we have proved that $d_{(C_G(Q),e_Q)}^{(u,e_R)}\circ \Phi_{\hat{M}_Q}=\omega_{\hat{V}_Q}((u,u))\bar{\Phi}_{\hat{M}_R}\circ d_{(C_H(Q),f_Q)}^{(u,f_R)}$; see (\ref{equation:compatibility R}). Multiplying $\varepsilon_Q$ on each side, we have $d_{(C_G(Q),e_Q)}^{(u,e_R)}\circ \varepsilon_Q\Phi_{\hat{M}_Q}=\varepsilon_Q\omega_{\hat{V}_Q}((u,u))\bar{\Phi}_{\hat{M}_R}\circ d_{(C_H(Q),f_Q)}^{(u,f_R)}$. By Proposition \ref{prop:compativility of omega} and the equality (\ref{equation:aaa}), we have 
		$$\varepsilon_Q\omega_{\hat{V}_Q}((u,u))=\omega_V((u_1,u_1))\omega_{\hat{V}_1}((u_2,u_2))\cdots\omega_{\hat{V}_{n-1}}((u_n,u_n))\omega_{\hat{V}_Q}((u,u))=\varepsilon_R.$$ 
		Summarising the above, the family $(\varepsilon_Q\cdot\Phi_{\hat{M}_Q})_{Q\subseteq P}$ of perfect isometries forms an isotypy in \cite[Definition 9.5.1]{Lin18b}, completing the proof.    $\hfill\square$
		
	}
\end{void}

\section{Proof of Proposition \ref{theo:counterexample}}\label{section:Proof for Q8}




We start to prove Proposition \ref{theo:counterexample}. Assume we are in the context of Proposition \ref{theo:counterexample}. So we can continue to use the notation in Notation \ref{notation:3.1} and Scetion \ref{section:On generalised decomposition maps}. Write
$$P=Q_8=\{\pm1,\pm \alpha,\pm \beta,\pm \gamma~|~(-1)^2=1,~\alpha^2=\beta^2=\gamma^2=\alpha\beta\gamma=-1\}.$$
An easy calculation shows that for any subgroup $Q$ of $P$ with $|Q|=4$, we have $C_P(Q)=Q$. Hence every subgroup of $P$ is fully $\F$-centralised. Assume that $k\otimes_\O V$ is a $3$-dimensional endotrivial $kP$-module. By \cite[Proposition 7.3.12]{Lin18b}, there are four choices of $V$ up to isomorphism. The irreducible ordinary character table of $P=Q_8$ is as follows:
$${\begin{array}{*{20}{c|ccccccc}}
				& \{1\} & \{-1\} & \{\alpha,-\alpha\} & \{\beta,-\beta\} & \{\gamma,-\gamma\}\\
			\hline 
			\chi_1& 1 & 1 & 1 & 1 & 1   \\
			\chi_2 & 1 & 1 & 1 & -1 &-1 \\
			\chi_3 & 1 & 1 & -1 & 1 &-1 \\
			\chi_4  & 1 & 1 & -1 & -1 & 1 \\
			\chi_5 & 2 & -2 &0 & 0 &0 \\
	\end{array}}$$
Since $V$ is an endotrivial $\O P$-module, we have $V\otimes_\O V^*\cong \O \oplus \O P$ as $\O P$-modules. Comparing the character values of the element $-1\in P$ on both sides we see that $(\rho_V(-1))^2=1$. By an easy calculation, we see that the four possibilities of $\rho_V$ are $\chi_1+\chi_5$, $\chi_2+\chi_5$, $\chi_3+\chi_5$ and $\chi_4+\chi_5$.
We may assume that $\rho_V=\chi_1+\chi_5$, and we can similarly prove Proposition \ref{theo:counterexample} for other three cases. Since the values of $\rho_V$ are in $\Z$, by Theorem \ref{theo:main}, $\Phi_M$ extends to an almost isotypy between $\O Gb$ and $\O Hc$. 

The set of all subgroups of $P$ is $\{\langle1\rangle,\langle -1\rangle, \langle \alpha\rangle,\langle \beta\rangle,\langle \gamma\rangle,P\}$. Since $V$ is an endotrivial $\O P$-module, by \cite[Lemma 3 (iii)]{Bilandadv}, we see that for any nontrivial subgroup $Q$ of $P$, $M_Q$ is a trivial source $kC_G(Q)e_Q$-$kC_H(Q)f_Q$-bimodules. Since $Q$ is fully $\F$-centralised, both $kC_G(Q)e_Q$ and $kC_H(Q)f_Q$ have $C_P(Q)$ as a defect group. By Proposition \ref{prop:vertex of slashed modules}, $(\Delta C_P(Q),e_{C_P(Q)}\otimes f_{C_P(Q)}^\circ)$ is a vertex subpair of $M_Q$.

Arguing by contradiction, suppose that there exists a family of signs $(\varepsilon_Q)_{\{1\neq Q\subseteq P\}}$ such that $\Phi_1:=\Phi_M$ together with $(\Phi_Q:=\varepsilon_Q\Phi_{\hat{M}_Q})_{\{1\neq Q\subseteq P\}}$ is an isotypy in \cite[Definition 9.5.1]{Lin18b}.

For any $1\neq Q\subseteq P$, since $(\Delta C_P(Q),e_{C_P(Q)}\otimes f_{C_P(Q)}^\circ)$ is a vertex subpair of $M_Q$, it is also a vertex subpair of $\hat{M}_Q$; see Remark \ref{remark:on source triples}. Let $\hat{V}_Q$ be a source of $\hat{M}_Q$ with respect to the vertex subpair $(\Delta C_P(Q),e_{C_P(Q)}\otimes f_{C_P(Q)}^\circ)$. Since $k\otimes_\O\hat{V}_Q$ is a source of $M_Q$, the $\O$-rank of $\hat{V}_Q$ is $1$. 
Note that the $\Z$-linear isomorphism $\bar{\Phi}_{\hat{M}_Q}:\Z\IBr_K(C_H(Q),\hat{f}_Q)\to \Z\IBr_K( C_G(Q),\hat{e}_Q)$ depends only on $M_Q$ - it does not depend on the choice of $\hat{V}_Q$ and $\hat{M}_Q$. Since $\rho_V(\alpha)=\rho_V(\beta)=\rho_V(\gamma)=1$ (recall that $\rho_V=\chi_1+\chi_5$) and $\rho_V(-1)=-1$,  by Notation \ref{notation:endopermutation modules} we have
\begin{equation}\label{equation:omega in section 8}
 \omega_V(\alpha)=\omega_V(\beta)=\omega_V(\gamma)=1~~{\rm and}~~\omega_V(-1)=-1.
 \end{equation}
  For any $u\in P$, since $\langle u\rangle$ is fully $\F$-centralised, by Proposition \ref{prop:bijections between simple modules of local blocks},
  $\bar{\Phi}_{\hat{M}_{\langle u\rangle}}(\varphi_{\alpha'})=\varphi_\alpha$ for any pair $(\alpha,\alpha')\in \mathcal{LP}(u,f_u)\times\mathcal{LP}(u,f_u)$ such that $\alpha\cap i\O Gi$ and $\alpha'\cap j\O Hj$ correspond under the bijection (\ref{equation:bijection of local pointed groups}).  Using the three equalities in \ref{void:decomposition numbers}, it is straightforward to check that
  \begin{equation*}
  	d_{(G,b)}^{(u,e_u)}(\Phi_M(\chi'))=\omega_V(u)\bar{\Phi}_{\hat{M}_{\langle u\rangle}}(d_{{(H,c)}}^{(u,f_u)}(\chi'))
  \end{equation*}
  for any $\chi'\in \Irr_K(kC_H(u),f_u)$.
  Hence we have
  \begin{equation}\label{equation:decomposition maps in section 8}
d_{(G,b)}^{(u,e_u)}\circ \Phi_M=\omega_V(u)\bar{\Phi}_{\hat{M}_{\langle u\rangle}}\circ d_{(H,c)}^{(u,f_u)}.
\end{equation}
Now (\ref{equation:omega in section 8}) and (\ref{equation:decomposition maps in section 8}) imply that $\varepsilon_{\langle \alpha\rangle}=\varepsilon_{\langle \beta\rangle}=\varepsilon_{\langle \gamma\rangle}=1$ and $\varepsilon_{\langle -1\rangle}=-1$. In other words, we have 
\begin{equation}\label{equation:8.1}
\Phi_{\langle \alpha\rangle}=\Phi_{\hat{M}_{\langle \alpha\rangle}},~~\Phi_{\langle \beta\rangle}=\Phi_{\hat{M}_{\langle \beta\rangle}},~~\Phi_{\langle \gamma\rangle}=\Phi_{\hat{M}_{\langle u\rangle}}~~{\rm and}~~\Phi_{\langle -1\rangle}=-\Phi_{\hat{M}_{\langle -1\rangle}}.
\end{equation}
By the transitivity of slash functors \cite[Lemma 22 (i)]{Biland}, $M_{\langle u\rangle}$ is a $(\Delta\langle u\rangle,e_u\otimes f_u^\circ)$-slashed module attached to $\hat{M}_{\langle -1\rangle}$ over $C_G(u)\times C_H(u)$. By Proposition \ref{prop:local source idempotents}, the bimodule $\hat{M}_{\langle-1\rangle}$ and the bimodule $M$ satisfy some similar conditions. Thus we can apply the equality (\ref{equation:decomposition maps in section 8}) to the $\O C_G(-1)\hat{e}_{\langle -1\rangle}$-$\O C_H(-1)\hat{f}_{\langle -1\rangle}$-bimodule $\hat{M}_{\langle -1\rangle}$ instead of the $\O Gb$-$\O Hc$-bimodule $M$, and we obtain
	\begin{equation}\label{equation:counterexample}
	d_{(C_G(-1),\hat{e}_{-1})}^{(u,e_u)}\circ \Phi_{\hat{M}_{\langle -1\rangle}}=\omega_{\hat{V}_{\langle -1\rangle}}((u,u))\bar{\Phi}_{\hat{M}_{\langle u\rangle}}\circ d_{(C_H(-1),\hat{f}_{-1})}^{(u,f_u)}.
	\end{equation}
Now consider the $\O \Delta C_P(-1)$-source $\hat{V}_{\langle -1\rangle}$ of $\hat{M}_{\langle-1\rangle}$. Since $C_P(-1)=P$, $\hat{V}_{-1}$ is an $\O \Delta P$-module. Since $P/[P,P]$ is a Klein four group and ${\rm rk}_\O(\hat{V}_{\langle -1\rangle})=1$, the $\O \Delta P$-module $\hat{V}_{\langle -1\rangle}$ is defined by a group homomorphism $P/[P,P]\to \O^\times$.
Now it is easy to see that there exists $u_0\in \{\alpha,\beta,\gamma\}$ such that $\omega_{\hat{V}_{\langle -1\rangle}}((u_0,u_0))=1$. Combining this with (\ref{equation:8.1}) and (\ref{equation:counterexample}), we have
	$$	d_{(C_G(-1),\hat{e}_{-1})}^{(u_0,e_{u_0})}\circ \Phi_{\langle -1\rangle}=-\bar{\Phi}_{\langle u_0\rangle}\circ d_{(C_H(-1),\hat{f}_{-1})}^{(u_0,f_{u_0})},$$
which contradicts to our assumption that  $\Phi_1:=\Phi_M$ together with $(\Phi_Q:=\varepsilon_Q\Phi_{\hat{M}_Q})_{\{1\neq Q\subseteq P\}}$ is an isotypy in the sense of \cite[Defnition 9.5.1]{Lin18b}.   $\hfill\square$


\bigskip\noindent\textbf{Acknowledgements.}\quad Many propositions and their proofs were obtained under Professor Markus Linckelmann's guidance, such as Proposition \ref{theo:counterexample}, Proposition \ref{prop:iMj}, Lemma \ref{lem:lift-uniqueness}, and so on. The author is very grateful to him for a lot of very useful discussions, and to City, University of London for its hospitality during the research for this paper in the spring of 2024. The author is also very grateful to Professor Yuanyang Zhou for bringing the author into the topic of \cite{HZ} in 2018.

\end{document}